\documentclass[preprint,11pt]{elsarticle}
\usepackage{geometry}
\newcommand{\beq}{\begin{equation}}
\newcommand{\eeq}{\end{equation}}
\newcommand{\bsq}{\begin{subequations}}
\newcommand{\esq}{\end{subequations}}
\newcommand{\bq}{\begin{eqnarray}}
\newcommand{\eq}{\end{eqnarray}}
\newcommand{\bqn}{\begin{eqnarray*}}
\newcommand{\eqn}{\end{eqnarray*}}
\usepackage{bbm}
\usepackage{enumerate}
\usepackage{enumitem} 

\usepackage{mathptmx}       
\DeclareMathAlphabet{\mathcal}{OMS}{cmsy}{m}{n}
\usepackage{helvet}         
\usepackage{courier}        
\usepackage{type1cm}        
\usepackage{makeidx}         
\usepackage{graphicx}        
\usepackage{url} 
\usepackage{multicol}        
\usepackage[bottom]{footmisc}
\usepackage{ifthen}
\usepackage{subfigure}
\usepackage[usenames,dvipsnames]{color}

\makeindex             

\usepackage{graphics} 
\usepackage{graphicx} 
\usepackage{epsfig} 
\usepackage{epstopdf}
\usepackage{mathptmx} 
\usepackage{times} 
\usepackage[cmex10]{amsmath}

\usepackage{mathtools}  
\usepackage{amssymb}  
\usepackage{amsthm}

\usepackage{amsfonts}
\usepackage{verbatim}
\usepackage{comment}
\usepackage{algorithmic}
\usepackage{multirow}
\usepackage{stfloats}
\usepackage{supertabular}
\usepackage{longtable}

\DeclareGraphicsExtensions{.pdf,.png,.jpg,.eps,.ps}

\usepackage{arydshln}
\usepackage{multicol}
\usepackage{colortbl}
\theoremstyle{definition}

\newtheorem{lemma}{Lemma}

\newtheorem{proposition}{Proposition}

\theoremstyle{definition}

\usepackage{comment}
\usepackage{ifthen}
\newboolean{showcomments}
\setboolean{showcomments}{true}
\newcommand{\ychen}[1]{\ifthenelse{\boolean{showcomments}}
        { \textcolor{red}{YC: #1}}}

\usepackage{amsmath}
\allowdisplaybreaks[4]

\usepackage[ruled]{algorithm2e}

\begin{document}

\begin{frontmatter}

\title{An Online Algorithm for Combined Computing Workload and Energy Coordination Within A Regional Data Center Cluster}
                
\author[Address1]{Shihan Huang}
\author[Address1]{Dongxiang Yan}             
\author[Address1]{Yue~Chen\corref{mycorrespondingauthor}}
\cortext[mycorrespondingauthor]{Corresponding author}
\ead{yuechen@mae.cuhk.edu.hk}

\address[Address1]{Department of Mechanical and Automation Engineering, The Chinese University of Hong Kong, Hong Kong, China.}

\begin{abstract}
Regional data center clusters have flourished in recent years to serve customers in a major city with low latency. The optimal coordination of data centers in a regional cluster has become a pressing issue because of its rising energy consumption. In this paper, a Lyapunov optimization-based online algorithm is developed for the combined computing workload and energy coordination of data centers in a regional cluster. The proposed online algorithm is prediction-free and easy to implement. We prove that the workload queues and battery energy level will be within their physical limits, though their related time-coupling constraints are not considered explicitly in the proposed algorithm. The previous online algorithms do not have such a guarantee. A theoretical upper bound on the optimality gap between the online and offline results is derived to provide a performance guarantee for the proposed algorithm. To enable distributed implementation, an accelerated ADMM algorithm is developed with iteration truncation and follow-up well-designed adjustments, whereby a nearly optimal solution is attained with much enhanced computational efficiency. Case studies show the effectiveness of the proposed method and its advantages over the existing methods.
\end{abstract}                


\begin{keyword}
data center \sep Lyapunov optimization \sep energy sharing \sep distributed coordination \sep online algorithm
\end{keyword}
                
\end{frontmatter}
        
\section*{NOMENCLATURE}
\subsection*{Abbreviations}
\begin{tabular}{ll}
DC & Data center. \\
PV & Photovoltaic.\\
MPC & Model predictive control. \\
P2P & Peer-to-peer. \\
ADMM & Alternating direction method of multipliers. \\
\end{tabular}

\subsection*{Indices and Sets}
\begin{tabular}{ll}
$i$ & Index of back ends. \\
$j,k$ & Index of front ends. \\
$\mathcal{I}$ & Set of all back ends. \\
$\Omega_i^B$ & Set of the front ends linked to back end $i$. \\
$\Omega_j^F$ & Set of the back ends linked to front end $j$. \\
\end{tabular}

\subsection*{Parameters}
\begin{tabular}{ll}
$A_j$ & Workload demands arrived at front end $j$. \\
$\underline{B}_i$ & Lower bound of battery energy level at back end $i$. \\
$\overline{B}_i$ & Upper bound of battery energy level at back end $i$. \\
$C_i$ & Upper bound of charging rate of the battery at back end $i$. \\
$D_i$ & Upper bound of discharging rate of the battery at back end $i$. \\
$E_i$ & Upper bound of the workloads back end $i$ processes. \\
$M_{ij}$ & Upper bound of the workloads transferred from front end $j$ to back end $i$. \\
$p^\text{b}_i$ & Electricity buying price from the main grid. \\
$p^\text{s}_i$ & Electricity selling price to the main grid. \\
$Q_i^{B}$ & Upper bound of the length of workload queue at back end $i$. \\
$Q_j^{F}$ & Upper bound of the length of workload queue at front end $j$. \\
$r_i$ & Constant in the definition of $l_i$. \\
$\overline{U}_{ik}$ & Upper bound of energy traded from back end $i$ to back end $k$. \\
$z_i$ & Energy of PV generation at back end $i$. \\
$\alpha_{ij}$ & Coefficient of bandwidth cost between back end $i$ and front end $j$. \\
$\beta_i$ & Coefficient of battery energy loss at back end $i$. \\
$\delta_i$ & Constant in the definition of $l_i$. \\
$\eta_\text{c}$ & Charge efficiency of battery. \\
$\eta_\text{d}$ & Discharge efficiency of battery. \\
$\gamma_j$ & Income from processing unit workload at front end $j$. \\ 
$\theta_j$ & Constant in the definition of $h_j^F$. \\
$\varphi_i$ & Constant in the definition of $h_i^B$. \\
\end{tabular}
    
\subsection*{Decision Variables}
\begin{tabular}{ll}
$a_j$ & Workloads accepted by front end $j$. \\
$b_i$ & Battery energy level at back end $i$. \\
$c_i$ & Charging rate of the battery at back end $i$. \\
$d_i$ & Discharging rate of the battery at back end $i$. \\
$e_i$ & Workloads back end $i$ processes. \\
$h_i^B$ & Virtual queue of workloads at back end $i$. \\
$h_j^F$ & Virtual queue of workloads at front end $j$. \\
$l_i$ & Virtual queue of battery energy at back end $i$. \\
$m_{ij}$ & Workloads transferred from front end $j$ to back end $i$. \\
$q_i^B$ & Length of workload queue at back end $i$. \\
$q_j^F$ & Length of workload queue at front end $j$. \\
$u_{ik}$ & Energy traded from back end $i$ to back end $k$. \\
$x_i$ & Electricity bought from the main grid at back end $i$. \\
$y_i$ & Electricity sold to the main grid at back end $i$. \\
\end{tabular}
        
\section{Introduction}
Recent decades have witnessed the boom of data centers (DCs) due to the information explosion. The installed base and the storage capacity of DCs globally were increased by 30\% and 26 times, respectively, during 2010-2018 \cite{miller2020sustainability}.
Regional DC clusters have become an inevitable trend in order to serve customers in major cities with low-latency \cite{regional}.
DCs are expected to be the largest energy consumers in 2030, contributing 
to 8\% of the global energy consumption \cite{WANG2022107926}. 
 To reduce the potential carbon emissions coming along, renewable generators such as photovoltaic (PV) panels have been installed in many DCs \cite{ALOBAIDI2021107231}. Most DCs are also equipped with energy storage to maintain a reliable supply of electricity. However, the supply-demand mismatch, caused by the significant spatial and temporal variation of computing workloads and renewable generations, cannot be fully offset even by battery storage. 
To reduce the consequent efficiency loss, the combined computing workload and energy coordination of DCs is necessary. Being clustered in a city level, regional DCs can not only share computing tasks with each other but also share energy to reduce their overall operation cost, which is the focus of this paper.

The coordination of DCs has captured great attention in recent years. 
To deal with the underlying uncertainties, stochastic optimization was used to minimize the expected operation costs of DCs based on a number of scenarios sampled from real data portfolio. 
References \cite{Ding} and \cite{WangP} focused on the day-ahead planning and real-time operation of DCs, respectively. 
The economic dispatch of power system with DCs was studied in \cite{NIU2021106358}, whereby the two-stage stochastic problem was solved by Benders decomposition. 
Robust optimization is another method of addressing uncertainties, by minimizing the worst-case DC operation cost \cite{Chen, Jawad}.   
Despite the fruitful research, the references above adopted offline models, which assume complete information over the whole time horizon. However, the accurate data of future electricity prices, renewable generations, and computing workloads are hardly available in practice. For actual implementation, the coordinated operation of DCs needs to be based on the up-to-date information without future predictions. Therefore, an online coordination method is necessary.

As a classic online optimization technique, model predictive control (MPC) has been applied in dynamic provisioning and load balancing in DCs with workloads assumed to be known in a look-ahead window \cite{Fang}.
The trade-off between the residential and migration costs of load balancing was handled by a randomized algorithm based on MPC technique in \cite{Mansouri}. 
An MPC-based control strategy was developed for DCs to participate in demand response programs \cite{cupelli2018data}. The MPC method is easy to implement, but it still requires certain predictions about the future. A large window is usually necessary to ensure optimality. Meanwhile, it may suffer from heavy computational burdens.

Lyapunov optimization is another online optimization technique that can be applied without prediction \cite{Guo}.
An event-triggered mechanism was proposed to control a regional integrated energy system by Lyapunov optimization in \cite{WANG2021106451}, which reduced the decision-making frequency of the system.
Lyapunov optimization was used in dynamic voltage and frequency scaling to promote the use of renewable energy, wherein an improvement in the service delay and deadline misses was shown \cite{Karimiafshar}.
Novel virtual queues of workloads were designed in \cite{Sun}, by which the worst-case service delay was bounded.
However, the previous research rarely considered the limitation on the backlog of workloads, which is determined by the storage capacity of DC. This is partly due to difficulty in applying Lyapunov optimization when the backlog limitation is considered, which this paper aims to overcome. 
Moreover, the works above did not consider the energy sharing among different DCs. In particular, the DCs in a regional cluster are located close to each other. By selling surplus energy to other DCs with energy shortage, the supply-demand mismatch in the region can be reduced and the overall operation efficiency of the regional cluster can be improved \cite{SONG2022108289}.

To facilitate energy sharing, a distributed implementation framework is needed for the purpose of privacy protection and the reduction of computational burdens. Several distributed algorithms have been developed for energy sharing among DCs \cite{Yu2}, charging stations \cite{yan2022distributed}, microgrids \cite{xu2020peer}, etc. One of the most widely-used algorithms is the alternating direction method of multipliers (ADMM) algorithm. However, a key deficiency of ADMM is its slow convergence, especially when the operational time scale is small. The convergence performance of ADMM was improved in some literature, e.g. by over-relaxation \cite{ZHANG2020106094} or by parallelization \cite{RAJAEI2021107126}. Reference \cite{ullah2021peer} tried to accelerate the ADMM algorithm by using a predictor-corrector technique. Unfortunately, these approaches all have their own shortcomings. For example, it may be tricky to tune the parameters of over-relaxation method for optimality.
When dealing with complex DC systems with fast-changing renewable generations, prediction and correction may still be too time-consuming.

This paper aims to develop a combined computing workload and energy coordination mechanism for DCs in a regional cluster, established in an online manner and able to adapt to the fast-changing environment. The computing workloads are properly allocated to the DCs, and the DCs in a regional cluster can share energy with each other.
Our main contributions are two-fold:

1)	\textbf{A prediction-free online algorithm.} A Lyapunov optimization-based algorithm is developed for the online combined computing workload and energy coordination of DCs. The proposed algorithm is prediction-free, meaning that the coordination decisions in each period are made based solely on up-to-date information, which is more practical. Distinct from the other existing online algorithms, a novel virtual queue is designed and used in the proposed algorithm. This allows our algorithm to ensure that the workload backlogs and battery energy level remain within their physical ranges, even though those range constraints are relaxed in the online algorithm. To the best of our knowledge, the current works do not have such a guarantee. An upper bound on the optimality gap between the online and offline algorithms is also derived.

2) \textbf{An accelerated ADMM algorithm.} For the distributed implementation of the online coordination algorithm above, an accelerated ADMM algorithm is developed. To be specific, the ADMM algorithm is sped up by stopping the iteration when a certain threshold is reached, followed by a set of well-designed adjustments to balance the computing workload and energy. The proposed algorithm can achieve a nearly optimal solution with much enhanced computational efficiency. Simulation results show that the computational time is reduced by 61\% compared to the traditional ADMM algorithm. This is critical to the operation of DC system where parameters such as renewable generations are likely to change very quickly.

The rest of this paper is organized as follows. Section \ref{sec:model} introduces the problem formulation. An online algorithm based on Lyapunov optimization is proposed in Section \ref{sec:online}, and an accelerated ADMM-based distributed coordination mechanism is developed in Section \ref{sec:distributed}. Section \ref{sec:result} presents the simulation results. Finally, the conclusions are drawn in \ref{sec:conclu}.

\section{Mathematical Formulation} \label{sec:model}

\subsection{System Overview}

\begin{figure}[t]
  \centering
  \includegraphics[width=0.5\textwidth]{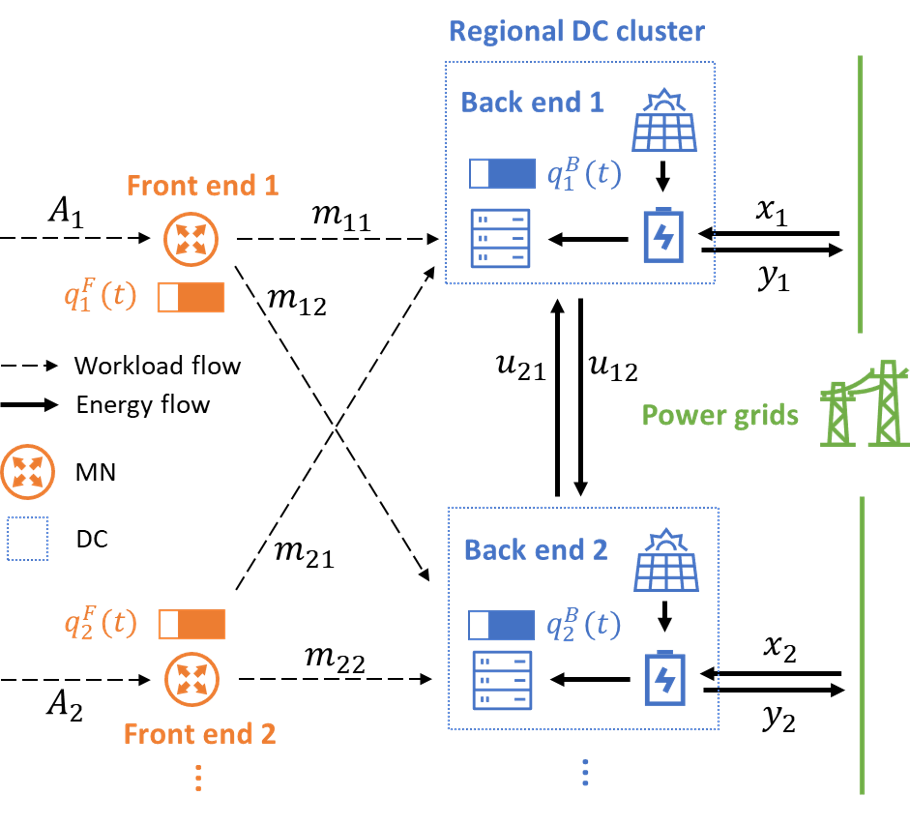}\\
  \caption{Illustration of the workload and energy flows among DCs.}\label{fig:sysConf}
\end{figure}

In this paper, we consider the coordination of a regional cluster of DCs, which are allowed to purchase electricity from the power grids and exchange energy with each other to serve the computing workloads. An illustration of the workload and energy flows among DCs is shown in Figure \ref{fig:sysConf}, consisting of the three elements below:
\begin{itemize}
    \item \emph{Front Ends} (mapping nodes, MNs), which accept computing workload demands and distribute them to different back ends for further processing.
    \item \emph{Back Ends} (DCs), which process workloads transferred from the front ends. Each back end is equipped with PV panels and battery energy storage. It can purchase electricity from the power grid and other DCs in the same regional cluster, use the power generated by its PV panels and energy storage to obtain the electricity needed to process the computing workloads. Surplus energy can be sold back to the power grid, used to charge the battery energy storage or shared with other DCs in the same regional cluster.    
    \item \emph{The Power Grid}, which provides electricity and buys back when necessary. The power grid operator also serves as a coordinator of the energy sharing among back ends.
\end{itemize}

As discussed above, due to the heterogeneity of renewable power generations and computing workloads, as well as the underlying uncertainties, online coordination of data centers is a must. In the following section, we start with an offline and centralized operation model for all DCs, then an online algorithm and a distributed coordination mechanism are developed in Sections \ref{sec:online} and \ref{sec:distributed}, respectively.



\subsection{Modeling of Data Center}
Suppose there are $I$ back ends and $J$ front ends, and the time interval investigated is divided into $T$ time slots. The data center related modelings are as follows.

\subsubsection{Workload Flow Related Modeling}
Denote the set of back ends linked to the front end $j$ and the set of front ends linked to the back end $i$ by $\Omega_j^F$ and $\Omega_i^B$, respectively.
For simplicity, we only consider non-critical requests, that is, part of the workloads arriving at the cluster can be rejected by the DCs, which will be either redirected to other clusters or ignored.
The workloads accepted by front ends will first be pushed into queues and then be transferred to other queues at back ends, waiting for processing.
All the workloads queues follow the first-in-first-out principle at front ends and back ends. The workload queues at the back end $i$ and the front end $j$ are denoted by $q_i^B(t)$ and $q_j^F(t)$, respectively. For $t = 1, 2, ..., T - 1,$ we have
\bsq \label{eq:workload}
\begin{align}
    q_j^F(t+1) = q_j^F(t) + a_j(t) - \sum \nolimits_{i \in \Omega_j^F} m_{ij}(t),\forall j, \label{cons10} \\
    q_i^B(t+1) = q_i^B(t) + \sum \nolimits_{j \in \Omega_i^B} m_{ij}(t) - e_i(t),\forall i, \label{cons11}
\end{align}
\esq
where $a_j(t)$ is the workloads accepted by the front end $j$ at the time slot $t$, $m_{ij}(t)$ is the workloads transferred from the front end $j$ to the back end $i$, $e_i(t)$ is the workloads that the back end $i$ process. The workload queues are initialized to zero, i.e., $q_j^F(1)=0,\forall j$, $q_i^B(1)=0,\forall i$. Due to the physical limitation of DCs, $a_j(t), m_{ij}(t)$ and $e_i(t)$ are bounded, i.e.,
\begin{gather}
    0 \le a_j(t) \le A_j(t), \forall j, \label{cons:a} \\
    0 \le m_{ij}(t) \le M_{ij}, \forall j, \forall i \in \Omega_j^F, \label{cons:m} \\
    0 \le e_i(t) \le E_i,\forall i, \label{cons:e}
\end{gather}
where $A_j(t)$ is the workload request at the front end $j$, $M_{ij}$ and $E_i$ are the upper bounds of $m_{ij}$ and $e_i$.

Similarly, the workload queues should be also bounded:
\bsq\label{eq:workloadbound}
\begin{align}
    0 \le q_j^F(t) \le Q_j^{F}, \forall j, \label{cons:qf} \\
    0 \le q_i^B(t) \le Q_i^{B}, \forall i, \label{cons:qb}
\end{align}
\esq
where $Q_j^F$ and $Q_{i}^B$ are the upper bounds of $q_j^F$ and $q_{i}^B$ respectively.

Transferring workload $m_{ij}(t)$ from the front end $j$ to the back end $i$ causes bandwidth cost, which is proportional to the amount $m_{ij}(t)$, i.e., 
\beq
    f^\text{tran}_{ij} (t) = \alpha_{ij} m_{ij}(t), \forall j, \forall i \in \Omega_j^F,
\eeq
where $\alpha_{ij}$ is the corresponding cost coefficient. The income of the front end $j$ is determined by the workloads it takes in the time slot $t$. Equivalently, we can use the disutility caused by rejecting part of the workloads for measurement:
\beq
    f_j^\text{work} (t) = \gamma_j (A_j(t) - a_j(t)), \forall j,
\eeq
where $\gamma_j$ is the disutility for not processing one unit of workload; $A_j(t)$ is the workload request arrived at the front end $j$, i.e., the maximum workloads front end $j$ can take. 

\subsubsection{Energy Flow Related Modeling}
In the time slot $t$, the back end $i$ buys electricity from the grids at a price $p^\text{b}_i(t)$ and sells back to the grids at a price $p^\text{s}_i(t)$, where $p^\text{b}_i(t) \ge p^\text{s}_i(t)$. Denote the amount of electricity bought and sold as $x_i(t)$ and $y_i(t)$ respectively, then the net cost of buying/selling electricity is:
\begin{gather}
    f_i^\text{grid} (t) = p^\text{b}_i(t) x_i(t) - p^\text{s}_i(t) y_i(t), ~ x_i(t) \ge 0, ~ y_i(t) \ge 0 \label{cons:x_y}, \forall i.
\end{gather}

Apart from the energy trading with the power grid, a back end can also charge or discharge its battery energy storage. Let $c_i(t)$ and $d_i(t)$ be the charging and discharging power of battery energy storage in back end $i$, respectively. Then,
\begin{align}
    0 \le c_i(t) \le C_i, ~ 0 \le d_i(t)\le D_i, ~ \forall i, \label{cons:chrg_dchg}
\end{align}
where $C_i$ and $D_i$ are the maximum charging and discharging power respectively. The energy loss caused by the charging and discharging of the battery energy storage is proportional to the charging and discharging power: 
\beq
    f_i^\text{batt} (t) = \beta_i (c_i(t) + d_i(t)), \forall i,
\eeq
where $\beta_i$ is the battery cost coefficient. 

The energy levels of battery in two consecutive time slots satisfy:
\begin{gather}
    b_i(t+1) = b_i(t) + \eta_\text{c} c_i(t) - \frac{1}{\eta_\text{d}} d_i(t), ~\forall i, ~ t = 1, 2, ..., T - 1, \label{cons8}
\end{gather}
where $\eta_\text{c}$ and $\eta_\text{d}$ are the charging and discharging efficiency respectively; $b_i(t)$ is the energy level of battery at the back end $i$ in the time slot $t$, bounded by:
\beq
    \underline{B}_i \le b_i(t) \le \overline{B}_i, ~ \forall i. \label{cons:batt_energy}
\eeq

Since all the DCs are located in a regional cluster, they can share energy with one another. Denote the energy sold by the back end $i$ to the back end $k$ by $u_{ik}(t)$, which satisfies the following constraints:
\bsq \label{eq:energy-share}
\begin{align}
      ~ & u_{ik}(t) + u_{ki}(t)= 0, ~ \forall k \in \mathcal{I} \backslash i, \label{cons:exchange} \\
    ~ & u_{ii}(t)=0,\forall i \in \mathcal{I}, \label{cons:exchange-2}\\
      ~&  -\overline{U}_{ik} \le u_{ik}(t) \le \overline{U}_{ik}, ~ \forall k \in \mathcal{I} \backslash i, \label{cons:u}
\end{align}
\esq
where $\mathcal{I} \backslash i$ denotes the set of all back ends except the back end $i$, $\overline{U}_{ik}$ is the maximum energy that can be shared. $u_{ik}>0$ means that the back end $i$ buys energy from the back end $k$. Otherwise, if $u_{ik}<0$, the back end $i$ sells to the back end $k$. Therefore, the shared energy satisfies the coupling constraint \eqref{cons:exchange}. Suppose the energy sharing price between back ends is $p^\text{t}$, with $p^\text{s}_i \le p^\text{t} \le p^\text{b}_i$ to encourage local energy exchange.

The power balance condition of the back end $i$ in the time slot $t$ is
\begin{align} 
    e_i(t) = ~ & x_i(t) - y_i(t) + d_i(t) - c_i(t) + z_i(t) + \sum_{\forall k \in \mathcal{I} \backslash i} u_{ik}(t), \forall i, \label{cons:back_end_engy}
\end{align}
where $z_i(t)$ is the PV generation in the back end $i$. 

The total costs of the overall system in the time slot $t$ is:
\begin{align}\label{eq:cost}
    f(t) = & \sum_{i = 1}^I \left(f_i^\text{grid} (t) + f_i^\text{batt} (t) \right) + \sum_{j = 1}^J \left(\sum_{i \in \Omega_j^F} f_{ij}^\text{tran} (t) + f_j^\text{work} (t) \right).
\end{align}
It is worth noting that the cost of energy sharing among back ends do not appear in \eqref{eq:cost} because the payment and income of back ends cancel each other according to \eqref{cons:exchange}.

\subsection{Offline and Centralized Operation of Data Centers}\label{subsec:form}
With the above constraints and objective, the offline and centralized operation of the overall system is formulated as:
\begin{align}
\textbf{P1: } \min ~ & \lim_{T \to \infty} \frac{1}{T}\sum_{t = 1}^T \mathbb{E}[f(t)]\label{obj}  \\
\mbox{s.t.} ~ & (\ref{eq:workload}) - (\ref{eq:workloadbound}), ~ (\ref{cons:x_y}), ~ (\ref{cons:chrg_dchg}) - (\ref{eq:energy-share}), \nonumber
\end{align}
where the objective \eqref{obj} aims to minimize the time-average expected total costs of the overall system.

Although \textbf{P1} is a linear program, it requires prior knowledge of electricity price $p^\text{b}_i (t)$, $p^\text{s}_i (t)$, workload demands $A_j(t)$, and PV generations $z_i(t)$ over the whole time horizon, which is difficult to obtain in practice. Therefore, an online counterpart that enables real-time decision-making is desired. In the following section, we propose an online algorithm for \textbf{P1}. 

Before we proceed, we make the following assumptions:



\noindent \textbf{A1}: $Q_j^F \ge A_{j,\max} + \sum_{i \in \Omega_j^F} M_{ij}, \forall j$,


\noindent \textbf{A2}: $Q_i^B \ge E_i+\sum_{j \in \Omega_i^B} M_{ij},\forall i$,

\noindent \textbf{A3}: $\overline{B}_i - \underline{B}_i \ge \eta_\text{c} C_i + \frac{1}{\eta_\text{d}}D_i,\forall i$,

\noindent \textbf{A4}: $P_{\max}^{\text{b}}\eta_{\text{c}}\eta_{\text{d}} < P_{\min}^{\text{s}} + \beta_i(1+\eta_\text{c}\eta_\text{d}),\forall i$,

\noindent where 
\bsq
\begin{align}
    A_{j,\max} &= \max\{A_j(t)~|~t = 1, 2, ..., T - 1\}, ~ \forall j \\
    E_{\max} &= \max \{E_i~|~ i \in \mathcal{I} \}, \\
    P_{\max}^\text{b} &= \max \{p^\text{b}_i (t)~|~ i \in \mathcal{I}, ~ t  = 1, 2, ..., T - 1\}, \\
    P_{\min}^\text{s} &= \min\{p^\text{s}_i(t)~|~ i \in \mathcal{I}, ~ t  = 1, 2, ..., T - 1\}.
\end{align}

\esq

\textbf{A1}--\textbf{A3} are mild assumptions that the capacity of front end queues, back end queues, and the battery are large enough. For example, if $\eta_\text{c}C_i>\frac{1}{\eta_\text{d}}D_i$, \textbf{A3} is satisfied if we have $\overline{B}_i-\underline{B}_i\ge 2\eta_\text{c} C_i$, which means that the battery can be charged at maximum power for two consecutive time slots. Similarly, if $\eta_\text{c}C_i<\frac{1}{\eta_\text{d}}D_i$, \textbf{A3} holds if the battery can be discharged at maximum power for two consecutive time slots. Either of the cases can be easily met. \textbf{A4} is natural which indicates that buying electricity from the grid, storing it in energy storage, and then selling it is not as good as not buying power at all, even if we sell at a higher price $P_{\max}^\text{b}$ and buy at a lower price $P_{\min}^\text{s}$. This assumption is to avoid unnecessary charging and discharging.

\section{Online Algorithm}\label{sec:online}



The main difficulty in developing an online counterpart lies in the time-coupling constraints (\ref{cons10}), (\ref{cons11}) and (\ref{cons8}) in \textbf{P1}.
In this section, we discuss how the time-coupling constraints can be equivalently removed by means of Lyapunov optimization. Then, by minimizing the upper bound of the drift-plus-penalty term, an online algorithm for solving \textbf{P1} can be derived.

\subsection{Problem Modification}
To adapt to the Lyapunov optimization framework, the time-coupling constraints (\ref{cons10}), (\ref{cons11}) and (\ref{cons8}) need to be reformulated in a time-average form. To be specific, summing \eqref{cons10} up over $t$ from 1 to $T - 1$ and divide it by $T$, we have
\begin{align}
    \frac{1}{T}q_j^F(T) = \frac{1}{T} q_j^F(1) + \frac{1}{T}\sum_{t=1}^T \left(a_j(t)-\sum_{i \in \Omega_j^F} m_{ij}(t)\right), \forall j
\end{align}
Since $q_j^F(1)$ and $q_j^F(T)$ are bounded, we have
\begin{align}
    \lim_{T \to \infty} \frac{1}{T}\sum_{t=1}^T \left(a_j(t)-\sum_{i \in \Omega_j^F} m_{ij}(t)\right) = 0,\forall j. \label{eq:time-average-1}
\end{align}
Similarly, constraint (\ref{cons11}) and (\ref{cons8}) can be transformed into
\begin{align}
    \lim_{T \to \infty} \frac{1}{T} \sum_{t=1}^T \left(\sum_{j \in \Omega_i^B} m_{ij}(t)-e_i(t)\right) = 0, \forall i, \label{eq:time-average-2} \\
    \lim_{T \to \infty} \frac{1}{T}\sum_{t=1}^T \left(\eta_\text{c} c_i(t)-\frac{1}{\eta_\text{d}}d_i(t)\right) = 0, \forall i. \label{eq:time-average-3}
\end{align}


Replacing the time-coupling constraints (\ref{cons10}), (\ref{cons11}) and (\ref{cons8}) by \eqref{eq:time-average-1}--\eqref{eq:time-average-3}, \textbf{P1} can be transformed into \textbf{P1}$'$:
\begin{align}
\textbf{P1}': \min ~ & \lim_{T \to \infty} \frac{1}{T}\sum_{t = 1}^T \mathbb{E}[f(t)]  \\
\mbox{s.t.} ~ & (\ref{cons:a})-(\ref{eq:workloadbound}), ~ (\ref{cons:x_y}), ~ \eqref{cons:chrg_dchg}, ~ \eqref{cons:batt_energy}-(\ref{eq:energy-share}), \nonumber \\
& \eqref{eq:time-average-1}-\eqref{eq:time-average-3}. \nonumber
\end{align}

As a matter of fact, \textbf{P1}$'$ is a relaxation of \textbf{P1} since it is easy to prove that any solution feasible to \textbf{P1} is also feasible to \textbf{P1}$'$. \textbf{P1}$'$ is an optimization with time-average objective and constraints, which fits the Lyapunov optimization framework.

\subsection{Lyapunov Optimization Based Method} 
\textbf{P1}$'$ is still time-coupled due to constraints \eqref{eq:time-average-1}--\eqref{eq:time-average-3}. The next step is to transform \textbf{P1}$'$ into a time-decoupled form by means of Lyapunov optimization. We construct three virtual queues and formulate a relaxation of \textbf{P1}$'$ which does not have time-coupling constraints by drift-plus-penalty method.

\subsubsection{Construct Virtual Queues}
The first step is to construct virtual queues for the front end workload dynamics, the back end workload dynamics, and the battery dynamics, denoted by $h_j^F(t)$, $h_i^B(t)$, and $l_i(t)$, respectively. The traditional method \cite{Sun} is to build virtual queues $h_j^F(t)$, $h_i^B(t)$ by modifying the workload queues in \eqref{cons10}-\eqref{cons11} slightly:
\bsq
\begin{align}
    h_j^F(t+1)= \max \left\{h_j^F(t)+a_j(t)-\sum \nolimits_{i \in \Omega_j^F} m_{ij}(t),0 \right\},\forall j, \label{eq:tra-queue-1}\\
    h_i^B(t+1)= \max \left\{h_i^B(t)+\sum \nolimits_{j \in \Omega_i^B} m_{ij}(t)-e_i(t), 0 \right\},\forall i.  \label{eq:tra-queue-2}
\end{align}
\esq
This method is straightforward and easy to implement. However, it cannot ensure that the corresponding workload queues are within their physical upper bounds $Q_j^F$, $Q_i^B$. To overcome this shortcoming, in this paper novel virtual queues are proposed as follows:
\bsq\label{vir_que}
\begin{align}
    h_j^F(t)=~ & q_j^F(t) - \theta_j, \forall j, \label{vir_que2}\\ 
    h_i^B(t)=~ & q_i^B(t) - \varphi_i, \forall i, \label{vir_que1} \\   
    l_i(t)=~ & b_i(t) - (\delta_i + r_i V), \forall i, \label{vir_que3}
\end{align}
\esq
where $\theta_j$, $\varphi_i$, $\delta_i$, $r_i$, $V$ are parameters to be determined later.

According to the definitions above, it is easy to prove that queues $h_j^F(t)$, $h_i^B(t)$, and $l_i(t)$ are mean rate stable. For example, for the queue $h_j^F(t)$, according to \eqref{vir_que2}, we have
\begin{align}
    h_j^F(t+1) = h_j^F(t) + a_j(t) - \sum_{i \in \Omega_j^F} m_{ij}(t),\forall j.
\end{align}
Hence,
\begin{align}
     \lim_{T \to \infty} \frac{h_j^F(T)}{T} =~ & \lim_{T \to \infty} \frac{h_j^F(1)}{T} +\frac{1}{T}\sum_{t=1}^T \left(a_j(t) - \sum_{i \in \Omega_j^F} m_{ij}(t)\right) = 0, ~ \forall j.
\end{align}
The last equation is due to \eqref{eq:time-average-1}.
Similarly, we have
\begin{align}
    \lim_{T \to \infty} \frac{h_i^B(T)}{T}=0,~ \lim_{T \to \infty} \frac{l_i(T)}{T}=0, \forall i.
\end{align}
which means that all these virtual queues grow slower than linearly over time.

\subsubsection{Obtain Lyapunov Function and Drift-Plus-Penalty}

Based on the virtual queues, we define the Lyapunov function as follows:
\begin{align}
    L(\boldsymbol{\Theta}(t)) = ~ & \frac{1}{2} \sum_{i = 1}^I \left(\left(h_i^B(t)\right)^2 + \left(l_i(t)\right)^2\right)  + \frac{1}{2} \sum_{j = 1}^J \left(h_j^F(t)\right)^2,
\end{align}
where $\boldsymbol{\Theta}(t) = (h_i^B(t),\forall i; h_j^F(t), \forall j; l_i(t), \forall i)$, representing the current state of the system, then Lyapunov drift from the time slot $t$ to $t + 1$ is defined as
\beq
    \Delta(\boldsymbol{\Theta}(t)) = \mathbb{E} \left[L(\boldsymbol{\Theta}(t + 1)) - L(\boldsymbol{\Theta}(t)) | \boldsymbol{\Theta}(t) \right].
\eeq

Then, the drift-plus-penalty term is given by
\begin{align}\label{eq:drift-plus-penalty}
\Delta(\boldsymbol{\Theta}(t))+V\mathbb{E}[f(t)|\boldsymbol{\Theta}(t)].
\end{align}
By minimizing \eqref{eq:drift-plus-penalty} instead of the original objective function \eqref{obj}, the time-average constraints \eqref{eq:time-average-1}--\eqref{eq:time-average-3} can be omitted as long as an appropriate $V$ is chosen, which will be discussed later. However, now the term $\Delta(\boldsymbol{\Theta}(t))$ in the objective becomes time-coupled. Hence, further transformation is needed.

\subsubsection{Minimizing the Upper Bound}
The idea is to derive a time-decoupled upper bound of \eqref{eq:drift-plus-penalty} and minimize it instead. According to (\ref{vir_que1}) -- (\ref{vir_que3}), we have
\begin{align}\label{ineq3}
    L(\boldsymbol{\Theta}(t+1)) - L(\boldsymbol{\Theta}(t))
    \le ~ & \sum_j h_j^F(t) \left(a_j(t) - \sum_{i \in \Omega_j^F} m_{ij}(t)\right) 
    + \sum_i h_i^B(t) \left(\sum_{j \in \Omega_i^B} m_{ij}(t)-e_i(t)\right) \nonumber\\
    ~ & + \sum_i l_i(t)\left(\eta_\text{c} c_i(t)-\frac{1}{\eta_\text{d}} d_i(t)\right) +N_1 +N_2+N_3 \nonumber \\
    := ~ & g(t) + N_1 + N_2 + N_3, 
\end{align}
where
\bsq
\begin{align}
    N_1 & = \frac{1}{2} \sum_j \max \left\{\max_t A_j^2 (t), \left(\sum_{i \in \Omega_j^F} M_{ij}\right)^2 \right\}, \\
    N_2 & = \frac{1}{2} \sum_i \max \left\{E_i^2, \left(\sum_{j \in \Omega_i^B} M_{ij} \right)^2 \right\}, \\
    N_3 & = \frac{1}{2} \sum_i \max \left\{\left(\eta_\text{c} C_i\right)^2, \left(\frac{D_i}{\eta_\text{d}}\right)^2 \right\}.
\end{align}
\esq

Ignoring the constant terms in the objective function (which will not change the optimal solution), the online optimization problem can be formulated as follows:
\begin{align}
\textbf{P2: } \min ~ & g(t) + V f(t) \label{obj_online} \\
\mbox{s.t.} ~ & (\ref{cons:a}) - (\ref{cons:e}), ~ (\ref{cons:x_y}), ~ \eqref{cons:chrg_dchg}, ~ (\ref{cons:back_end_engy}) - (\ref{eq:energy-share}). \nonumber
\end{align}
Note that constraints bounding workload queues \eqref{eq:workloadbound} and the energy level of battery \eqref{cons:batt_energy} have been removed, and the time-average constraints \eqref{eq:time-average-1}--\eqref{eq:time-average-3} are no longer needed because all of them are automatically satisfied by minimizing \eqref{obj_online} as proved in Proposition \ref{propo1} later. There is no time-coupling constraint in \textbf{P2}, which can be computed online.

\subsection{Performance Guarantee}
To ensure that the optimal solution of \textbf{P2} is feasible to \textbf{P1}, the parameter $V$ in \eqref{eq:drift-plus-penalty} and the parameters in \eqref{vir_que} must satisfy certain requirements as follows.

\begin{proposition} \label{propo1}
When assumptions A1--A4 hold, if $V$ and $r_i$ satisfies the following requirements:
\bsq \label{eq:requirement}
\begin{align}
    V \ge & \max_j \max_{i \in \Omega_j^F} \frac{\sum_{i \in \Omega_j^F} M_{ij}-\theta_j + \varphi_i}{\alpha_{ij}}, \label{v_req1} \\
    V \ge & \max_i \left\{\max_{j \in \Omega_i^B} \frac{Q_j^F - Q_i^B + \sum_{j \in \Omega_i^B} M_{ij} - \theta_j+\varphi_i}{\alpha_{ij}},
    \frac{E_i-\varphi_i}{P_{\min}^\text{s}}, ~
    \frac{\eta_\text{c} (\delta_i + \eta_\text{c} C_i - \overline{B}_i)} {P_{\min}^\text{s} + \beta_i - \eta_\text{c} r_i}, \frac{\frac{1}{\eta_\text{d}}(\underline{B}_i+\frac{1}{\eta_\text{d}}D_i-\delta_i)} {-P_{\max}^\text{b} + \beta_i + \frac{1}{\eta_\text{d}}r_i} \right\}, \label{v_req2} \\
    V \le & \min_j \left\{ \frac{Q_j^F - A_{j,\max} - \theta_j} {\gamma_j} \right\}, \label{v_req3}
\end{align}
\begin{flalign}
    & P_{\max}^\text{b} \eta_\text{d}-\beta_i\eta_\text{d} < r_i <  \frac{P_{\min}^\text{s} + \beta_i}{\eta_\text{c}}, \forall i, &
\end{flalign}
\esq
then the optimal solution obtained by \textbf{P2} satisfies the constraints (\ref{cons:qf}), (\ref{cons:qb}) and (\ref{cons:batt_energy}).
\end{proposition}

The proof of Proposition \ref{propo1} can be found in Appendix A. Due to \textbf{A4}, a proper $r_i$ is easy to find. To find the $V$ that satisfy \eqref{eq:requirement}, let
\bsq
\begin{align}
   \delta_i = ~ & \overline{B}_i - \eta_\text{c} C_i, \forall i, \\
   \theta_j= ~ & E_{\max} + \sum \nolimits_{i \in \Omega_j^F} M_{ij}, \forall j, \\
   \varphi_i=~ & E_i, \forall i.
\end{align}
\esq
Then for any $V > 0$, the requirement \eqref{v_req1} is always satisfied, and \eqref{v_req2} is reduced to
\begin{align}
    V \ge & \max_i \max_{j \in \Omega_i^B} \frac{1}{\alpha_{ij}} \left(Q_j^F - Q_i^B + \sum_{j \in \Omega_i^B} M_{ij}
    - \sum_{i \in \Omega_j^F} M_{ij} - (E_{\max} - E_i) \right). \label{v_req4}
\end{align}
Note that the right-hand side of \eqref{v_req4} will not be a large number, since the capacities of front ends $Q_j^F$ and back ends $Q_i^B$ are usually close and so is the bandwidth. Hence, it is not difficult to find a $V$ that satisfies \eqref{v_req3} and \eqref{v_req4}.

Apart from the satisfaction of time-coupling constraints, another issue we care about is the optimality gap between the offline problem \textbf{P1} and the online problem \textbf{P2}. 
Denote the values of $f(t)$ in the optimal solutions of \textbf{P1} and \textbf{P2} by $f^*(t)$ and $\hat{f}(t)$ respectively, the optimal value of \textbf{P1} by $F^*$, and let
\beq
    \hat{F} = \lim_{T \to \infty} \frac{1}{T} \sum_{t = 1}^T \mathbb{E} \left[\hat{f}(t)\right].
\eeq
The optimality gap can be bounded in the proposition below.
\begin{proposition} \label{propo2} The optimality gap is bounded by:
\begin{align} \label{ineq:gap}
    \hat{F} - F^* \le \frac{1}{V}\left(N_1 + N_2 + N_3\right).
\end{align}
\end{proposition}

The proof of Proposition \ref{propo2} can be found in Appendix B. The discussion above reveals that the parameter $V$ is critical to the performance of the proposed algorithm. According to Proposition \ref{propo2}, $V$ should be as large as possible to minimize the gap between $F^*$ and $\hat{F}$. Meanwhile, $V$ is bounded above by \eqref{v_req3}, otherwise the time-coupling constraints in \textbf{P1} will be violated.


\section{Distributed Implementation} \label{sec:distributed}
\textbf{P2} is still a centralized optimization problem, which may be impractical due to privacy concerns, communication and computational burdens. Therefore, a distributed algorithm is needed. In this section, an accelerated ADMM-based algorithm with iteration truncation is developed.

First, to distinguish the local variables of the three subsystems, we use $m_{ij}(t)$, $m'_{ij}(t)$ to denote the transferred workloads optimized by front ends and back ends respectively, and $u_{ik}(t)$, $u'_{ik}(t)$ to denote the exchanged energy optimized by back ends and the power grid respectively. Constraints (\ref{cons:exchange})-\eqref{cons:exchange-2} is then replaced by the following constraints:
\bsq
\begin{align}
    u'_{ik}(t) + u'_{ki}(t) & = 0, ~ \forall i, \forall k \in \mathcal{I} \backslash i, \label{cons:exchange_prime} \\
    u'_{ii}(t) & = 0, ~ \forall i. \label{cons:exchange_prime_1}
\end{align}
\esq
The two new variables should be also bounded:
\bsq
\begin{align}
    0 \le m'_{ij}(t) \le m_{ij}, ~ \forall j, \forall i \in \Omega_j^F, \label{cons:m_prime} \\
    \underline{U}_{ik} \le u'_{ik}(t) \le \overline{U}_{ik}, ~ \forall i, \forall k \in \mathcal{I} \backslash i. \label{cons:u_prime_1}
\end{align}
\esq
Coupling constraints are added to ensure that the optimal strategies provided by different subsystems are equal:
\bsq
\begin{align}
    & m'_{ij}(t) = m_{ij}(t), ~ \forall j,\forall i \in \Omega_j^F, \label{cons:m_prime_1} \\
    & u'_{ik}(t) = u_{ik}(t), ~ \forall i, \forall k \in \mathcal{I} \backslash i. \label{cons:u_prime_2}
\end{align}
\esq

In addition, $g(t)$ in \eqref{obj_online} also needs modification. Since the second term is optimized by back ends, $m_{ij}$ in this term is replaced by $m'_{ij}$. Denote the modified $g(t)$ by $g'(t)$:
\begin{align}
    g'(t) = & \sum \nolimits_j H_j^F(t) \left(a_j(t) - \sum \nolimits_{i \in \Omega_j^F} m_{ij}(t) \right)
    + \sum \nolimits_i H_i^B(t) \left(\sum \nolimits_{j \in \Omega_i^B} m'_{ij}(t) - e_i(t) \right) \nonumber \\
    & + \sum \nolimits_i B_i(t) \left(\eta_\text{c} c_i(t) - \frac{1}{\eta_\text{d}} d_i(t) \right).
\end{align}
Then the augmented Lagrangian function is given by
\begin{align} \label{obj_distr}
    \mathcal{L}(t) = ~ & g'(t) + Vf(t)
    + \frac{\rho}{2} \sum_j \sum_{i \in \Omega_j^F} \left(m_{ij}(t) -  m'_{ij}(t) + \lambda_{ij}(t) \right)^2
    + \frac{\rho}{2} \sum_i \sum_{k \in \mathcal{I} \backslash i} \left(u_{ik}(t) - u'_{ik}(t) + \mu_{ik}(t) \right)^2,
\end{align}
where $\lambda_{ij}(t)$, $\mu_{ik}(t)$ are the corresponding dual variables. Then \textbf{P2} can be equivalently transformed into:
\begin{align} \label{eq:P3}
   \textbf{P3: } \min ~ & \mathcal{L}(t) \\
   \mbox{s.t.} ~ & (\ref{cons:a}) - (\ref{cons:e}), ~ (\ref{cons:x_y}), ~ \eqref{cons:chrg_dchg}, ~ \eqref{cons:exchange-2}, ~ \eqref{cons:u}, ~ (\ref{cons:back_end_engy}), \nonumber \\
   & \eqref{cons:exchange_prime} - \eqref{cons:u_prime_2}. \nonumber
\end{align}


The traditional ADMM algorithm can be applied to solve \textbf{P2} in a distributed manner. In particular, the decision variables of front ends, back ends and the power grids, denoted by $X_{F,j}(t),\forall j$, $X_{B,i}(t),\forall i$ and $X_P(t)$ respectively, are 
\bsq
\begin{align}
    X_{F,j}(t) & = \{a_j(t), m_{ij}(t) ~|~ i \in \Omega_{j}^F\},\forall j, \\
    X_{B,i}(t) & = \{m'_{ij}(t), e_i(t), x_i(t), y_i(t), c_i(t), d_i(t),
    u_{ik}(t) ~|~ j \in \Omega_{i}^B, k \ne i\},\forall i, \\
    X_P(t) & = \{u'_{ik}(t) ~|~ k \ne i\},\forall i.
\end{align}
\esq
In the $n$-th iteration, the front end $j$ solves
\begin{align}\label{eq:front-end-problem}
    \min_{X_{F,j}(t)}~ & H_j^F(t)\left(a_j(t) - \sum_{i \in \Omega_j^F} m_{ij}(t) \right)+V\sum_{i \in \Omega_j^F} f_{ij}^\text{tran}(t) \\
    ~ &+Vf_j^\text{work}(t)+\frac{\rho}{2} \sum_{i \in \Omega_j^F} \left(m_{ij}(t) -  m'_{ij}(t) + \lambda_{ij}(t) \right)^2, \nonumber\\
    \mbox{s.t.}~ & \eqref{cons:a}, ~\eqref{cons:m}. \nonumber
\end{align}
The back end $i$ solves
\begin{align}\label{eq:back-end-problem}
   \min_{X_{B,i}(t)}~ &  H_i^B(t) \left(\sum_{j \in \Omega_i^B} m'_{ij}(t) - e_i(t) \right)
   + B_i(t) \left(\eta_\text{c} c_i(t) - \frac{1}{\eta_\text{d}} d_i(t) \right) + V (f_i^\text{grid}(t)+f_i^\text{batt}(t)) \\
   ~ & + \frac{\rho}{2} \sum_{j \in \Omega_i^B} \left(m_{ij}(t) -  m'_{ij}(t) + \lambda_{ij}(t) \right)^2
    + \sum_{k \in \mathcal{I} \backslash i} \left(p^\text{t}(t) u_{ik}(t)+ \frac{\rho}{2}\left(u_{ik}(t) - u'_{ik}(t) + \mu_{ik}(t) \right)^2\right), \nonumber\\
    \mbox{s.t.}~ & \eqref{cons:e}, ~\eqref{cons:x_y}, ~\eqref{cons:chrg_dchg}, ~\eqref{cons:u}, ~\eqref{cons:back_end_engy}, ~\eqref{cons:m_prime}. \nonumber
\end{align}
The power grid solves
\begin{align}\label{eq:grid-problem}
    \min_{X_P(t)} ~&\sum_i \sum_{k \in \mathcal{I} \backslash i} \left(u_{ik}(t) - u'_{ik}(t) + \mu_{ik}(t) \right)^2, \\
    \mbox{s.t.}~& \eqref{cons:exchange_prime}, ~\eqref{cons:exchange_prime_1}, ~\eqref{cons:u_prime_1}. \nonumber
\end{align}

However, the traditional ADMM may need considerable iterations to converge. The random parameters, e.g. electricity price and workloads, may have already changed before the algorithm is converged, rendering it less practical. A straightforward idea to lighten the computational burden is to break the iteration when it takes too much time to converge. However, this may lead to violation of the coupling constraints \eqref{cons:m_prime_1} and \eqref{cons:u_prime_2}. To tackle this issue, we propose the following adjustment method.

Suppose the iteration is truncated at the step $N$, we choose $(u'_{ik})^{(N)}$ as the optimal strategy and let
\beq\label{eq:48}
   {u}_{ik}^{(N)}(t) = (u'_{ik})^{(N)}(t).
\eeq
Then, constraint \eqref{cons:m_prime_1} is met but constraint (\ref{cons:back_end_engy}) may be violated due to the change of $u_{ik}^N(t)$. Denote by $\Delta e_i(t)$ the gap between the left-hand side and the right-hand side of (\ref{cons:back_end_engy}):
\begin{align} \label{eq:post_proc1}
    \Delta e_i(t) = & ~ x_i^{(N)}(t) - y_i^{(N)}(t) + d_i^{(N)}(t) - c_i^{(N)}(t)
    + z_i^{(N)}(t) + \sum_{\forall k \in \mathcal{I} \backslash i} (u'_{ik})^{(N)}(t) - e_i^{(N)}(t). 
\end{align}
This gap can be filled by adjusting the amount of energy bought from or sold to the grids. If $\Delta e_i(t) < 0$, we increase $x_i(t)$ by the same amount; otherwise we increase $y_i(t)$ by the same amount:
\bsq
\begin{align}
    \hat{x}_{i}(t) = x_{i}^{(N)}(t) - \min \{\Delta e_i(t), 0\}, \\
    \hat{y}_{i}(t) = y_{i}^{(N)}(t) + \max \{\Delta e_i(t), 0\}. \label{eq:post_proc2}
\end{align}
\esq

The treatment for $m_{ij}$ and $m'_{ij}$ is similar. Let
\begin{align}
   m_{ij}^{(N)}(t) = (m'_{ij})^{(N)}(t).
\end{align}
Then \eqref{cons:m_prime_1} is met but $q_j^F(t+1)$ may go beyond the bounds when updating the workload queue using \eqref{cons10}. To this end, ${a}^{(N)}_j(t)$ needs to be adjusted. If $q_j^F (t+1) < 0$, let
\bsq\label{eq:52}
\beq \label{eq:post_proc3}
    {a}^{(N)}_j(t) = \sum_{i \in \Omega_j^F} (m'_{ij})^{(N)}(t) - {q}_j^F(t).
\eeq
If $q_j^F (t+1) > Q_j^F$, let
\beq \label{eq:post_proc4}
    {a}^{(N)}_j(t) = \sum_{i \in \Omega_j^F} (m'_{ij})^{(N)}(t) - {q}_j^F(t) + Q_j^F.
\eeq
\esq
Then update $q_j^F(t+1)$ using the new $a^{(N)}_j(t)$.

The overall procedure is given in Algorithm \ref{algo1}. 

\begin{algorithm}[ht]
\caption{Accelerated ADMM Algorithm} \label{algo1}
\begin{algorithmic}[1]
\STATE Set the limit of iterations $N$, let $n=0$. 
\STATE Initialize $\lambda_{ij}^0(t)$ and $\mu_{ik}^0(t)$.
\REPEAT
\STATE For each front end $j$: Solve \eqref{eq:front-end-problem} to get $X_{F,j}^{(n + 1)}(t)$. 
\STATE For each back end $i$: Solve \eqref{eq:back-end-problem} to get $X_{B,i}^{(n + 1)}(t)$.
\STATE The power grid: Solve \eqref{eq:grid-problem} to get $X_P^{(n + 1)}(t)$.
\STATE{Update dual variables:
    \begin{align*}
        \lambda_{ij}^{(n + 1)}(t) & = \lambda_{ij}^{(n)}(t) + m_{ij}^{(n + 1)}(t) - (m'_{ij})^{(n + 1)}(t) \\
        \mu_{ij}^{(n + 1)}(t) & = \mu_{ij}^{(n)}(t) + u_{ik}^{(n + 1)}(t) - (u'_{ik})^{(n + 1)}(t)
    \end{align*}}
\UNTIL convergence \OR $n \ge N$
\STATE Adjustments according to \eqref{eq:48}--\eqref{eq:52}.
\end{algorithmic}
\end{algorithm}


\section{Case Studies} \label{sec:result}

In this section, the performance of the proposed algorithm is presented and compared with benchmark algorithms. Impact of various factors on the performance is also analyzed. Finally, the scalability of the proposed algorithm is examined.

\subsection{Simulation Setup} \label{subsec:setup}

The proposed algorithm is implemented on MATLAB 2022a and the simulations are performed on a desktop PC with an Intel i5-10505 CPU and 8 GB RAM. We first use a simple DC system with 2 front ends and 3 back ends for illustration. The parameters related to the batteries are $\overline{B}_i = 85$ MWh, $\underline{B}_i = 10$ MWh, $C = (4, 5, 3)$ kWh, $D = (6, 7.5, 4.5)$ kWh, $\eta_\text{c} = 0.9$, $\eta_\text{d} = 0.95$. The upper bounds of workload queues are $Q_i^B = Q_j^F = 200$ MWh. Real-world data are employed including the electricity prices of PJM \cite{pjm}, the workload traces of Google cluster \cite{google}, and solar radiation data of NREL \cite{solardata}. The simulation is conducted over a one-week period divided into 2016 time slots, where each time slot is 5 minutes.

To testify the effectiveness of the proposed online algorithm, we examine the performance of the following algorithms or models with the same parameters for comparison:
\begin{enumerate}[label=B\arabic*, font=\bfseries]
    \item Offline algorithm, i.e., solving \textbf{P1} directly assuming complete future information. \label{algo:offline}
    \item Greedy algorithm. We set a threshold $P_\text{th}$ of buying prices of electricity. The battery of back end $i$ is charged at the maximum power, i.e., $c_i(t) = C_i$, if the buying price at back end $i$ is lower than $P_\text{th}$. Then $f(t)$ is minimized in each time slot with the same constraints as \textbf{P1} for that specific time slot.  \label{algo:greedy}
    \item Proposed online algorithm on a model without energy sharing between back ends. This can be done by setting the lower and upper bounds of shared energy $\pm \bar U_{ik}$ to zero. \label{algo:no_exchange}
    \item Traditional online algorithm using virtual queues \eqref{eq:tra-queue-1}, \eqref{eq:tra-queue-2}, and \eqref{vir_que3}. \label{algo:unbounded}
\end{enumerate}

\subsection{Performance Comparisons}

\begin{figure}[htbp]
    \centering
    \begin{minipage}{0.4\textwidth}
        \includegraphics[width=6cm]{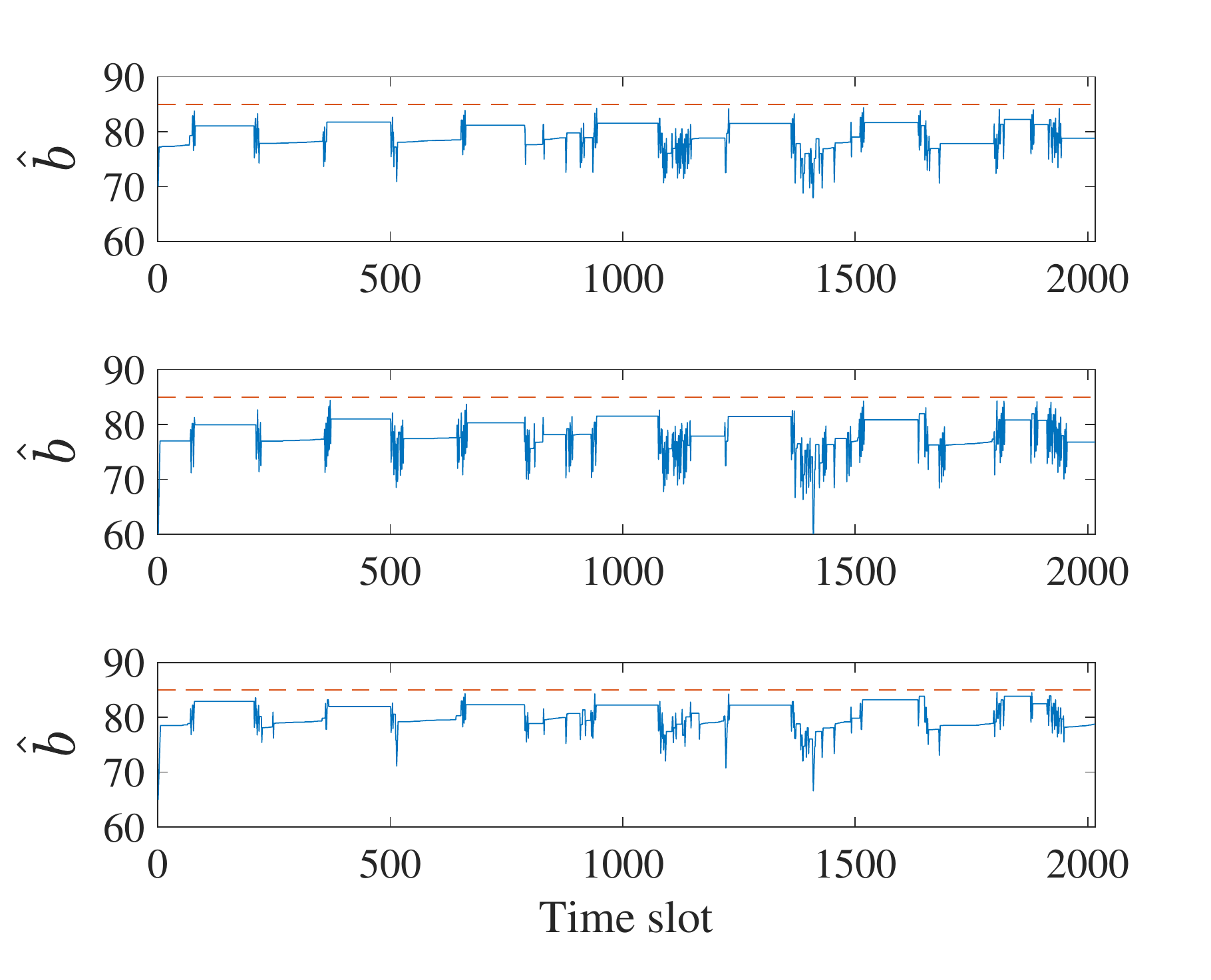}
    \end{minipage}
    \begin{minipage}{0.4\textwidth}
        \includegraphics[width=6cm]{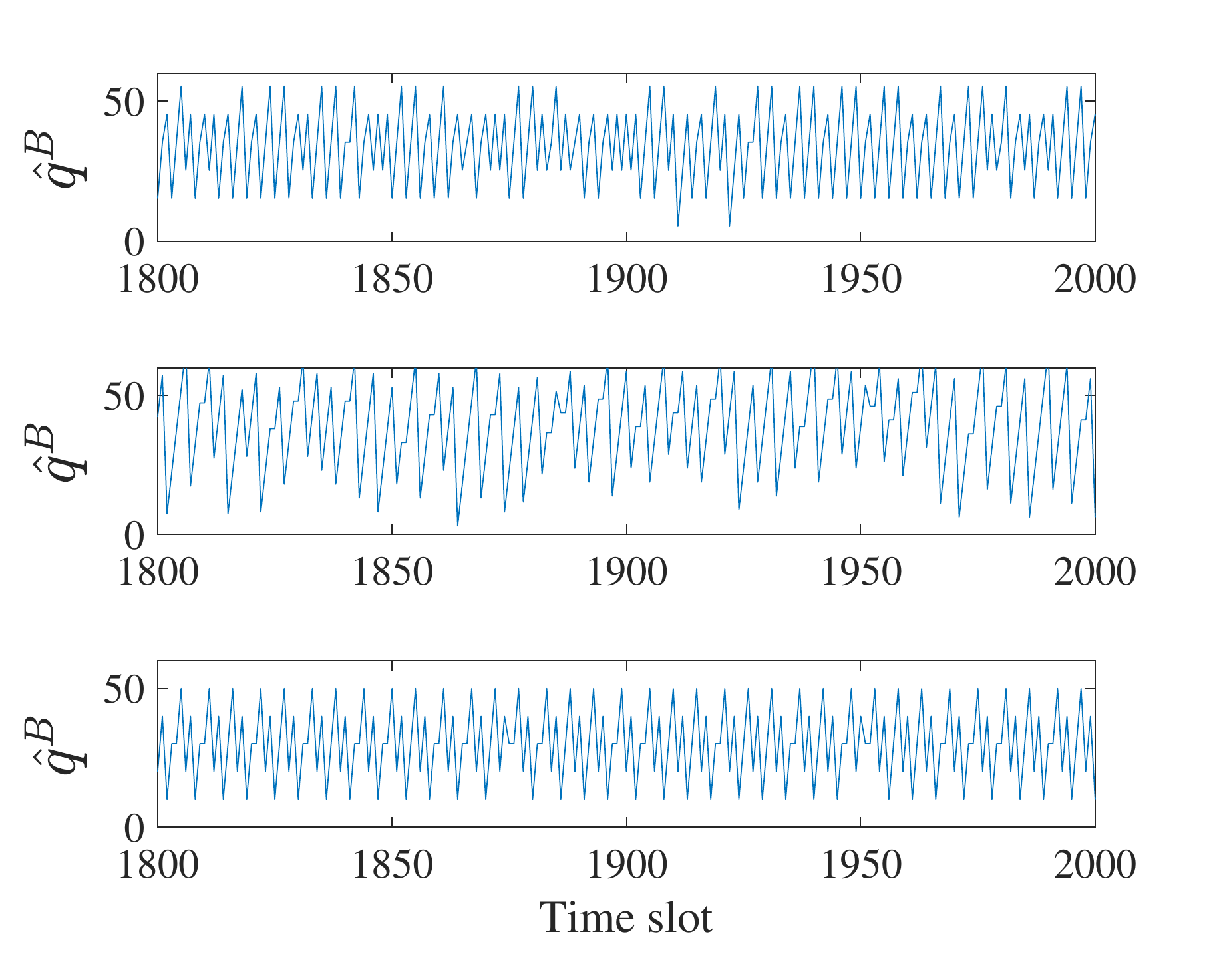}
    \end{minipage}
    \caption{Battery energy and workload queues of back ends. Blue line: battery energy trace; red dash line: upper bound of battery energy.} \label{fig:b_qb_hat}
\end{figure}

\begin{figure}[htbp]
    \centering
    \begin{minipage}{0.4\textwidth}
        \includegraphics[width=6cm]{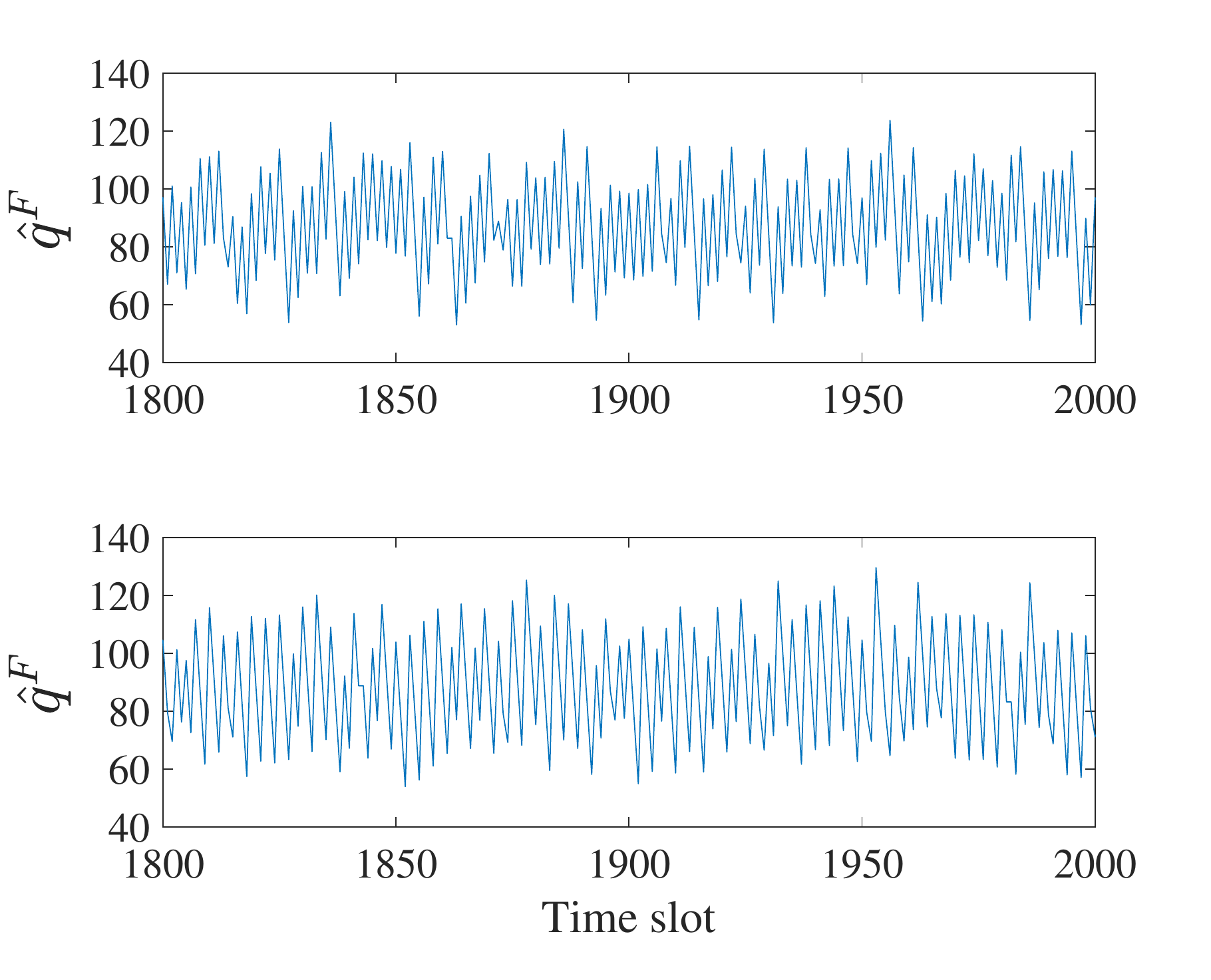}
    \end{minipage}
    \begin{minipage}{0.4\textwidth}
        \includegraphics[width=6cm]{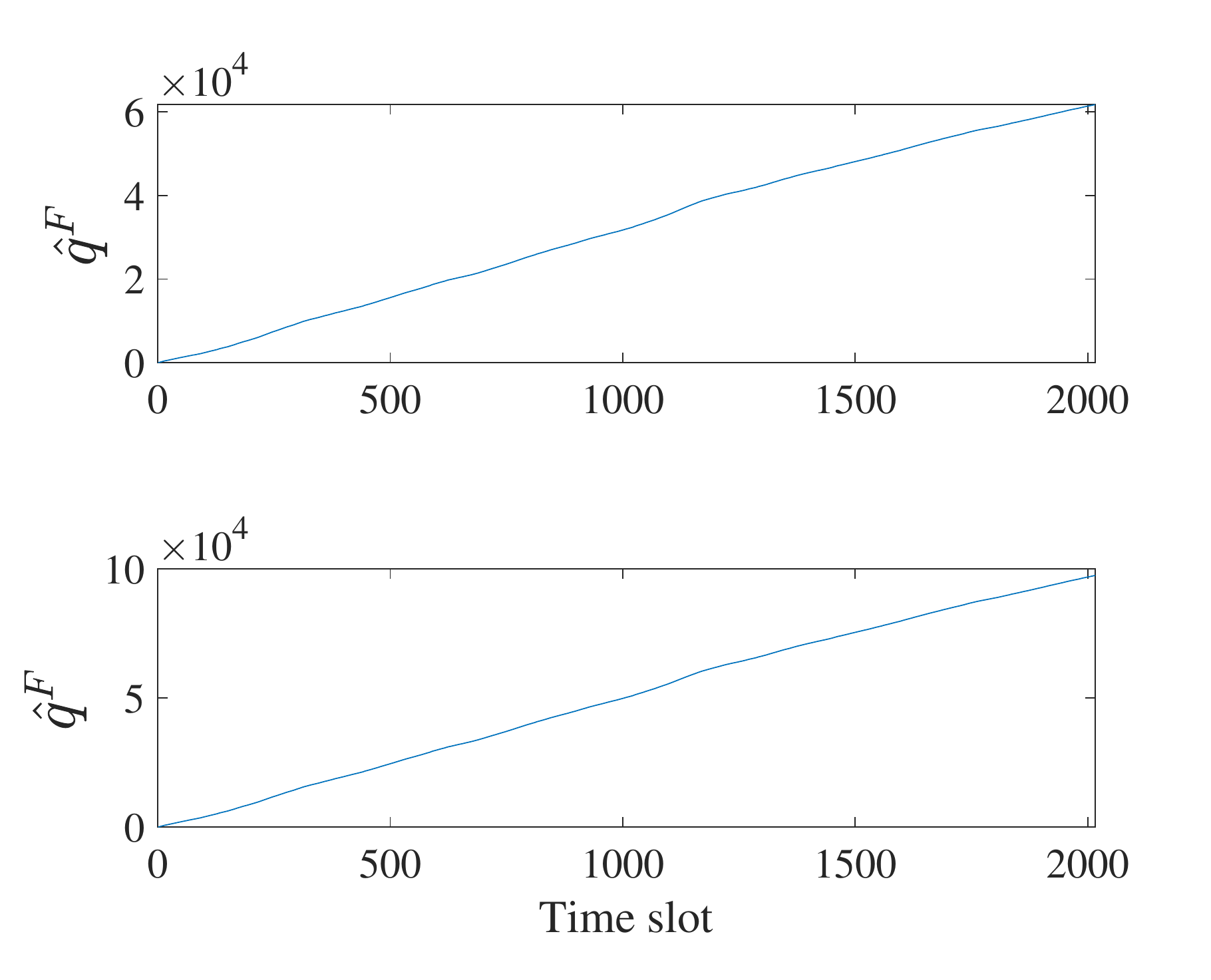}
    \end{minipage}
    \caption{The workload queues at front ends. Left: bounded workload queues by the proposed algorithm. Right: Unbounded workload queues by \ref{algo:unbounded}.} \label{fig:qf}
\end{figure}

Figure \ref{fig:b_qb_hat} shows the battery energy of the back ends using the proposed algorithm. It can be observed that constraint (\ref{cons:batt_energy}) is always met for all back ends. All traces are below the maximum battery energy (the red dash line), while the minimum battery energy is far below the plot window.
The workload queues at front ends and back ends are presented in Figure \ref{fig:b_qb_hat} (right) and Figure \ref{fig:qf} (left) respectively. Only part of the traces ($1800 \le t \le 2015$) is presented for the sake of clarity. Constraints (\ref{cons:qf}) and (\ref{cons:qb}) are all satisfied though they are not explicitly considered in \textbf{P2}, which justifies Proposition \ref{propo1}.

\begin{figure}[!htbp]
  \centering
  \includegraphics[width=0.45\textwidth]{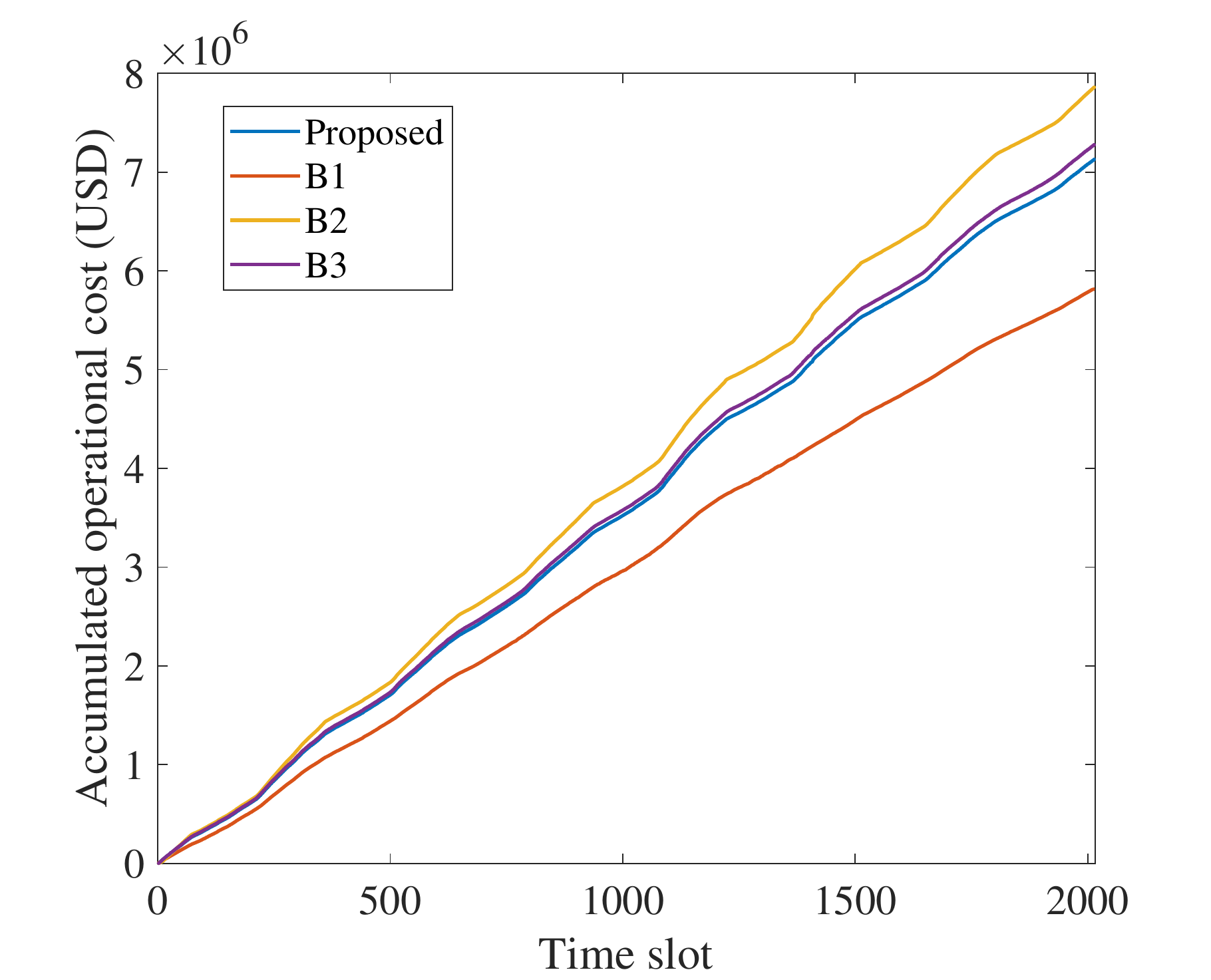} \\
  \caption{Accumulated operation costs under the four algorithms.} \label{fig:f_hat}
\end{figure}

Furthermore, we compare the performance of the proposed online algorithm with algorithms \ref{algo:offline}, \ref{algo:greedy} and \ref{algo:no_exchange} stated in Section \ref{subsec:setup}. The workload queues attained by \ref{algo:unbounded} is presented in Figure \ref{fig:qf} (right). In the results of B4, the queues at all front ends keep going up to far beyond their physical limitation, which reveals the necessity of the proposed algorithm that can ensure the satisfaction of physical limitation. The accumulated operation costs, i.e., the total costs from the beginning to the present time slot, of the overall DC system are shown in Figure \ref{fig:f_hat}. Compared to B2 and B3, the proposed algorithm reduces the total costs by 12\% and 2\% respectively. This shows the advantages of the proposed algorithm and that energy sharing among different data centers can improve the overall efficiency. Although the offline model \ref{algo:offline} has the lowest total costs, it is of little practical use due to the lack of future information. The time-average gap between the proposed algorithm and B1 is 655 USD at the end of the time horizon, less than the right-hand side of \eqref{ineq:gap}, which verifies Proposition \ref{propo2}. 
The details of the costs are listed in Table \ref{tab:costCmp}. 

\begin{table*}[!htbp]
\footnotesize
  \centering
    \caption{Cost comparison among B1, B2, B3 and proposed algorithm (Unit: Million USD).} \label{tab:costCmp}
   \begin{tabular}{@{}ccccccccc@{}}
     \hline\hline
                    & $f_\text{batt}$ & $f_\text{grid}$ & $f_\text{tran}$ & $f_\text{work}$ & $f$ & Relative value of $f$ \\
     \hline
     B1 & 0.002 & -0.026 & 1.178 & 4.660 & 5.813 & 100\% \\
     B2 & 0.441 & 2.105 & 2.217 & 3.102 & 7.864 & 135\% \\
    B3 & 0.045 & 1.556 & 1.509 & 4.171 & 7.281 & 125\% \\
     Proposed & 0.037 & 1.412 & 1.501 & 4.182 & 7.132 & 123\% \\
     \hline \hline
   \end{tabular}
\end{table*}

\subsection{Impact of Parameters}

\begin{figure}[htb]
    \centering
    \begin{minipage}{0.4\textwidth}
        \includegraphics[width=6cm]{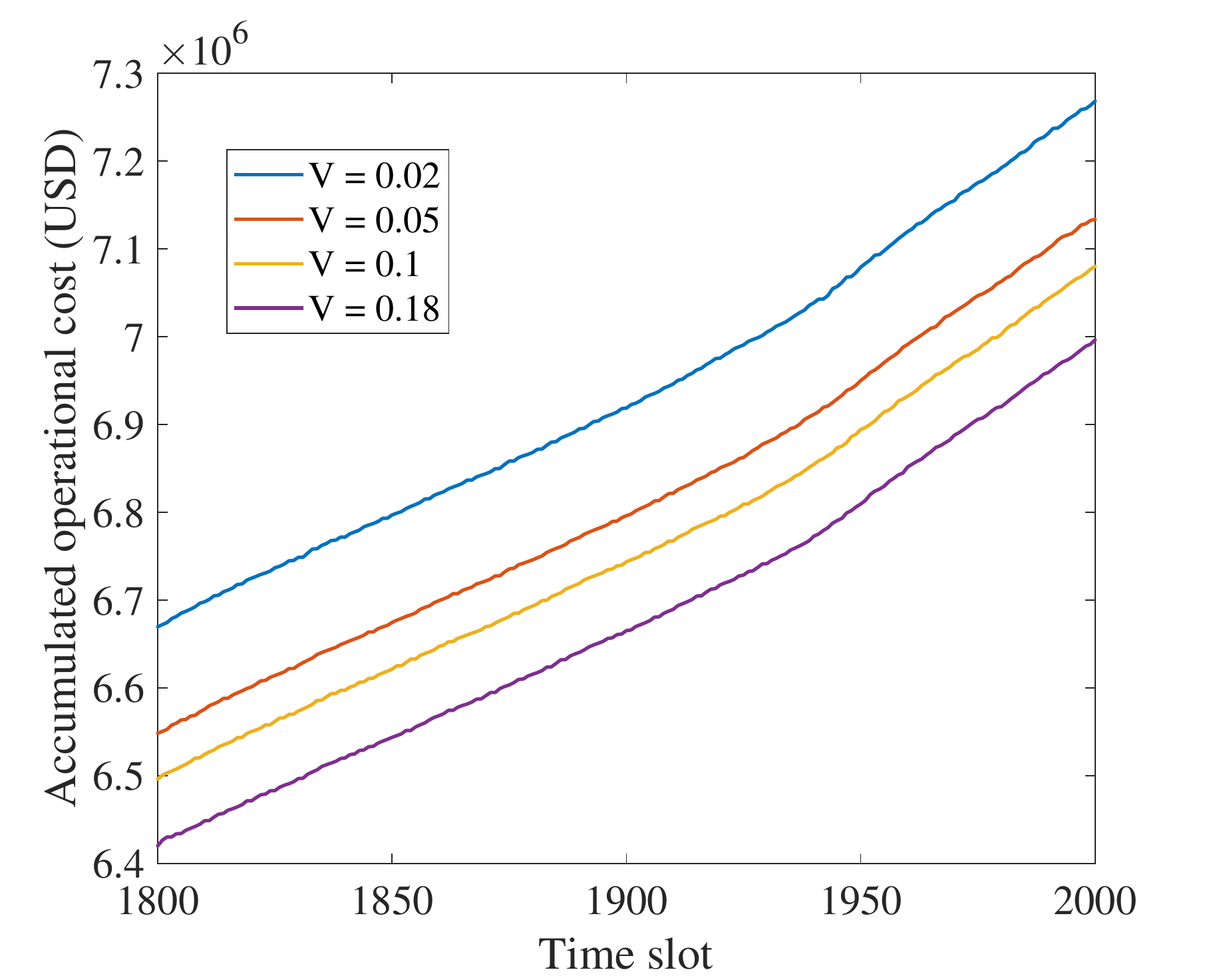}
    \end{minipage}
    \begin{minipage}{0.4\textwidth}
        \includegraphics[width=6cm]{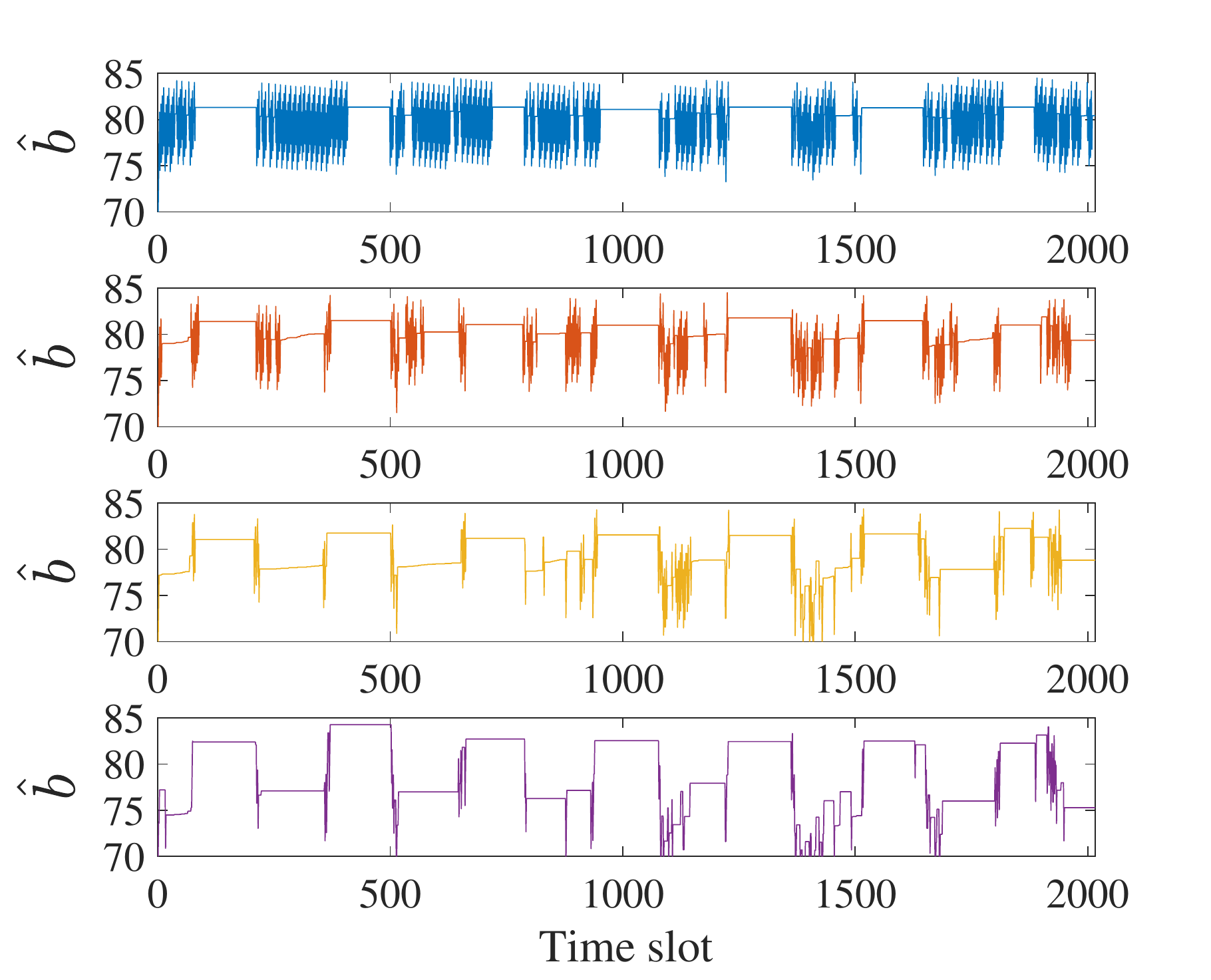}
    \end{minipage}
    \caption{The impact of $V$. Left: part of the traces of accumulated costs. Right: the cost of buying electricity at the end of time horizon.} \label{fig:compare_v}
\end{figure}

We further test the impact of $V$ in the objective function \eqref{obj_online} of \textbf{P2}. Let $V=0.02, 0.05, 0.1, 0.18$, respectively and record the accumulated overall operation costs in Figure \ref{fig:compare_v}. It can be observed that the costs are reduced as $V$ increases since it puts more emphasis on $f(t)$ in (\ref{obj_online}). The traces of battery energy provide some insights into this phenomenon: As shown in Figure \ref{fig:compare_v}, the charging and discharging of battery become more frequent as $V$ goes down. This is because the stability of virtual queues \eqref{vir_que1}--\eqref{vir_que3}, including the one related to battery energy, is prioritized since $g(t)$ dominates the objective function, keeping battery energy level in a smaller range. Consequently, the battery cost is increased from $2.35 \times 10^4$ USD when $V = 0.18$ to $1.57 \times 10^5$ USD when $V = 0.02$.

\begin{figure}[htb]
    \centering
    \begin{minipage}{0.4\textwidth}
        \includegraphics[width=6cm]{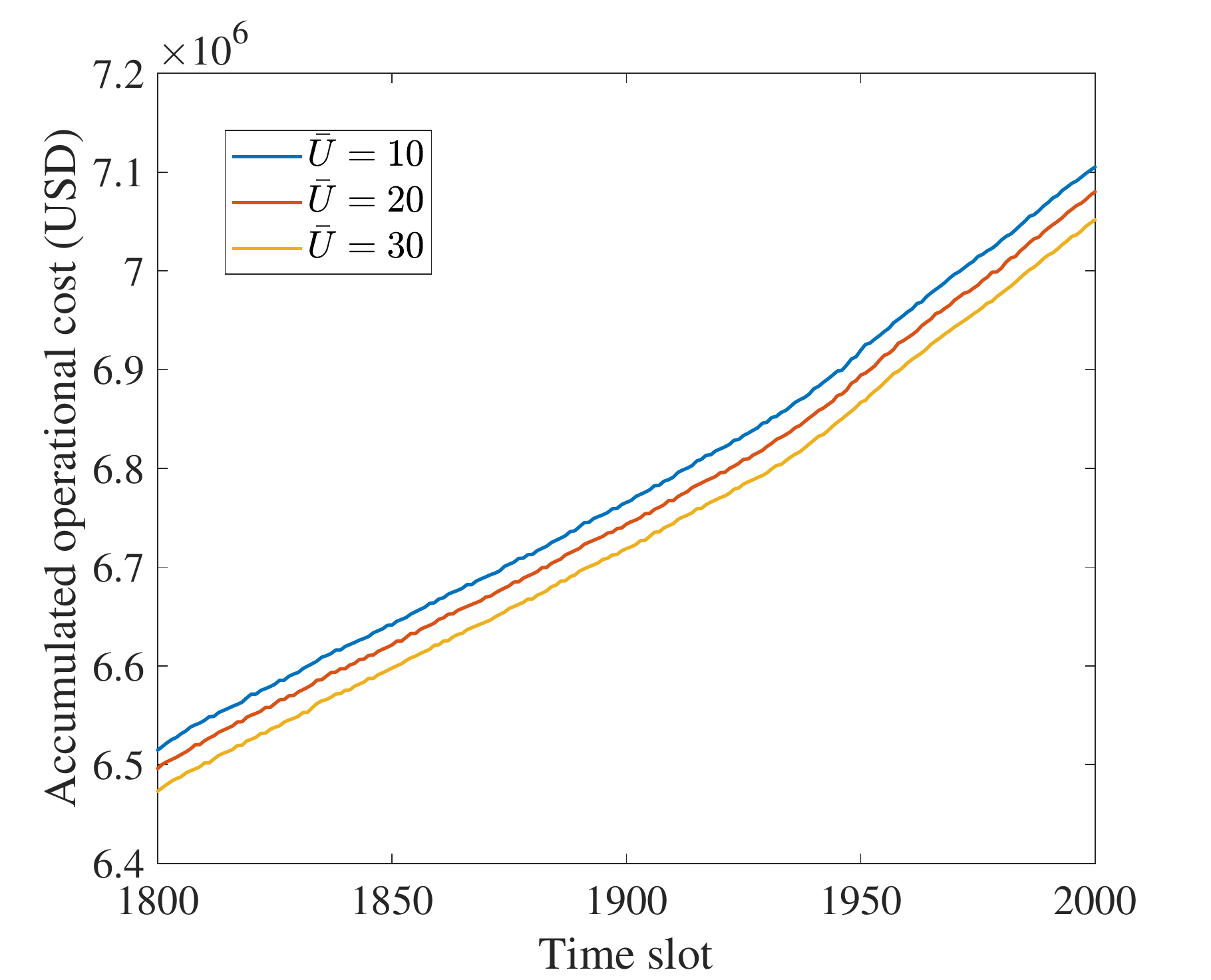}
    \end{minipage}
    \begin{minipage}{0.4\textwidth}
        \includegraphics[width=6cm]{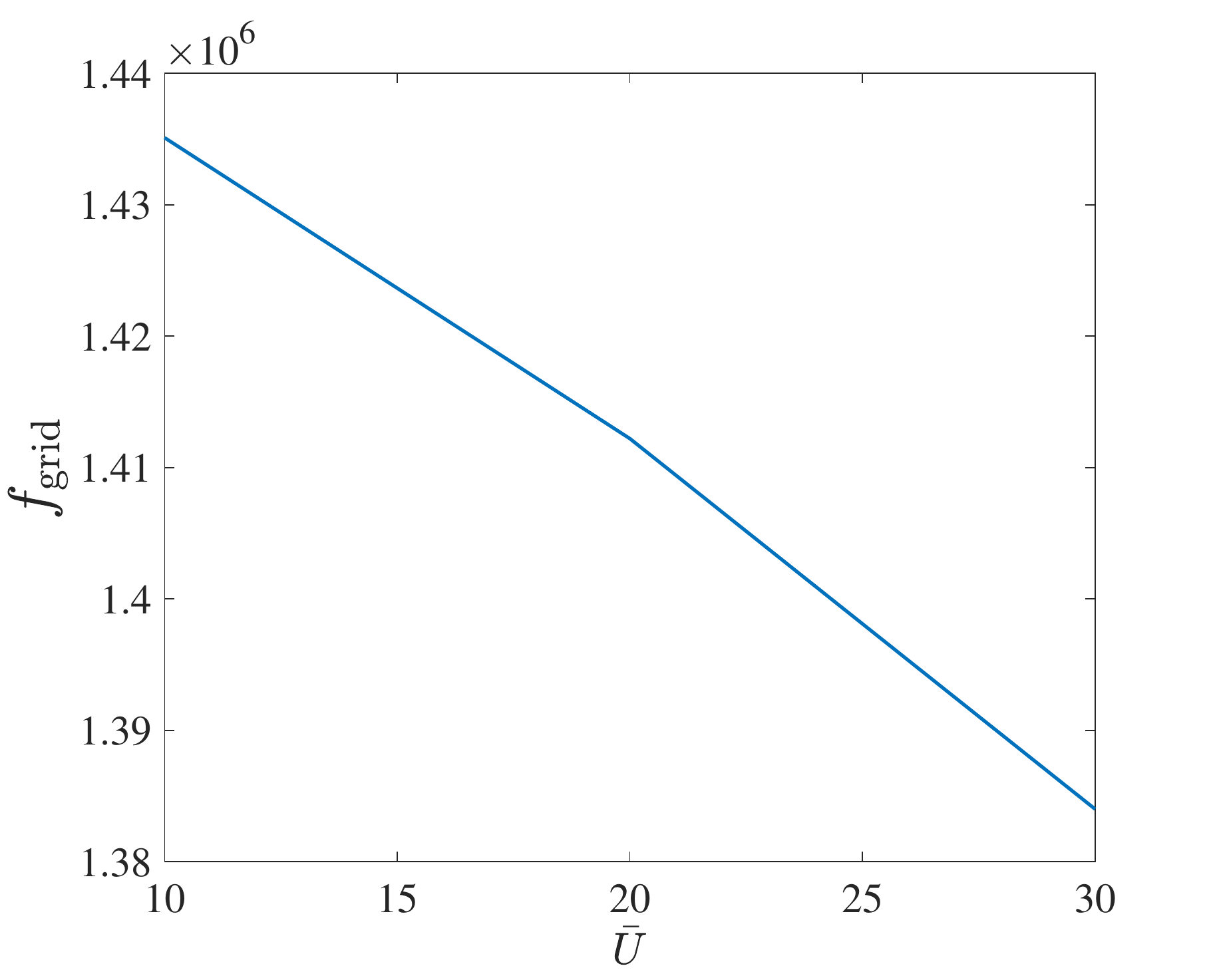}
    \end{minipage}
    \caption{The impact of $\bar U$. Left: part of the traces of accumulated costs. Right: the net cost of buying electricity at the end of time horizon.} \label{fig:compare_u}
\end{figure}

The role of P2P energy sharing on the total costs of DCs are also investigated. A drop in the costs can be observed in Figure \ref{fig:compare_u} as the upper bound of P2P trading $\bar U$ is increased, revealing the improvement in the efficiency of DC system. The reason is that: The DC short of energy may turn to other DCs instead of buying electricity from the main grids, since $p^\text{t} \le p^\text{b}_i$, the payment of the DC decreases. Similarly, since $p^\text{t} \ge p^\text{s}_i$, the DC with surplus energy can sell to other DCs instead of selling to the grid to earn a profit.
The relationship between the net cost of buying electricity $f_\text{grid}$ and $\bar U$ is presented in Figure \ref{fig:compare_u}.

\subsection{Convergence \& Scalability}

\begin{figure}[htbp]
  \centering
  \includegraphics[width=0.45\textwidth]{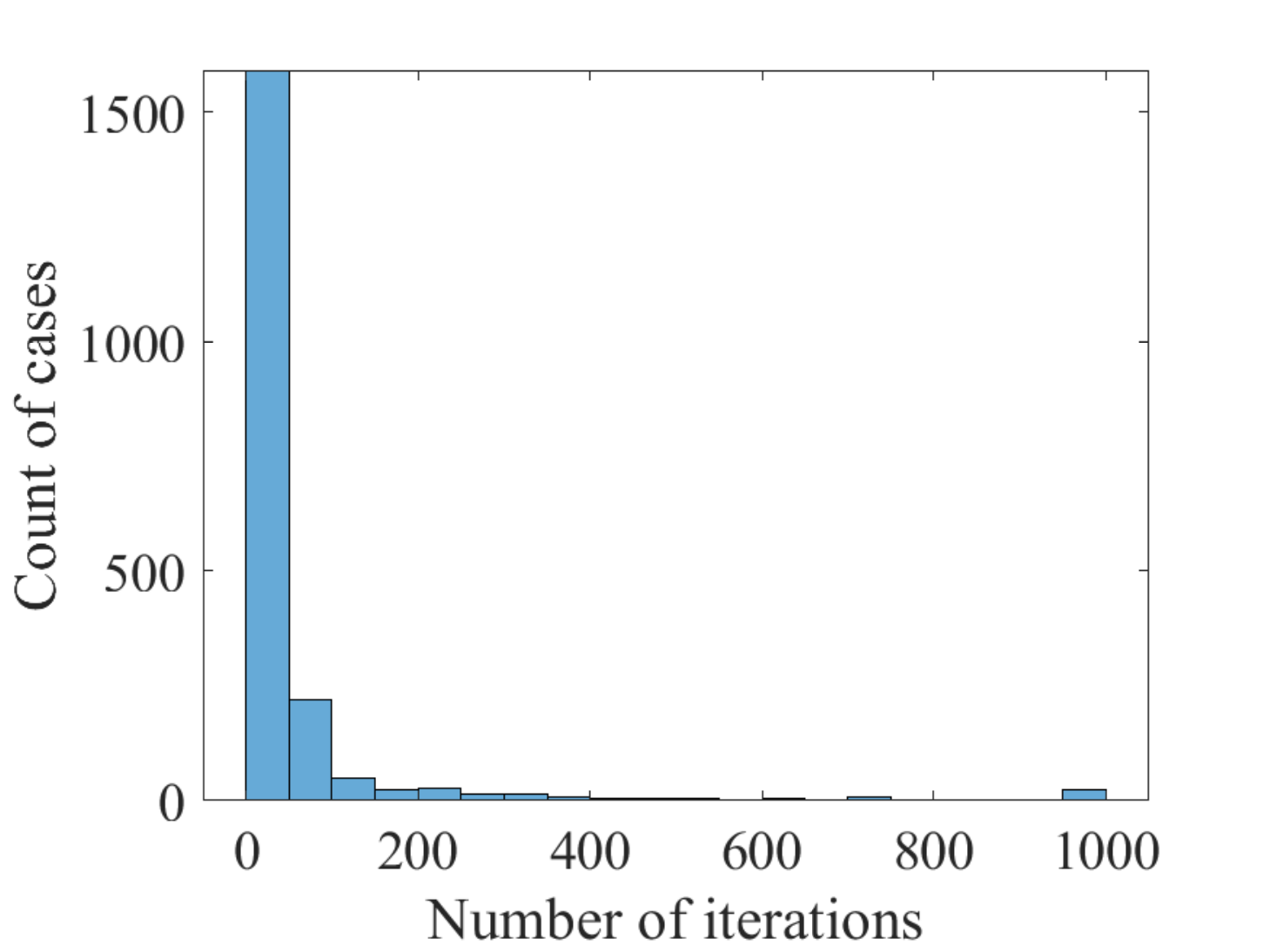} \\
  \caption{Histogram of the numbers of iterations.} \label{fig:histo}
\end{figure}

\begin{figure}[htbp]
  \centering
  \begin{minipage}{0.45\textwidth}
      \includegraphics[width=7.5cm]{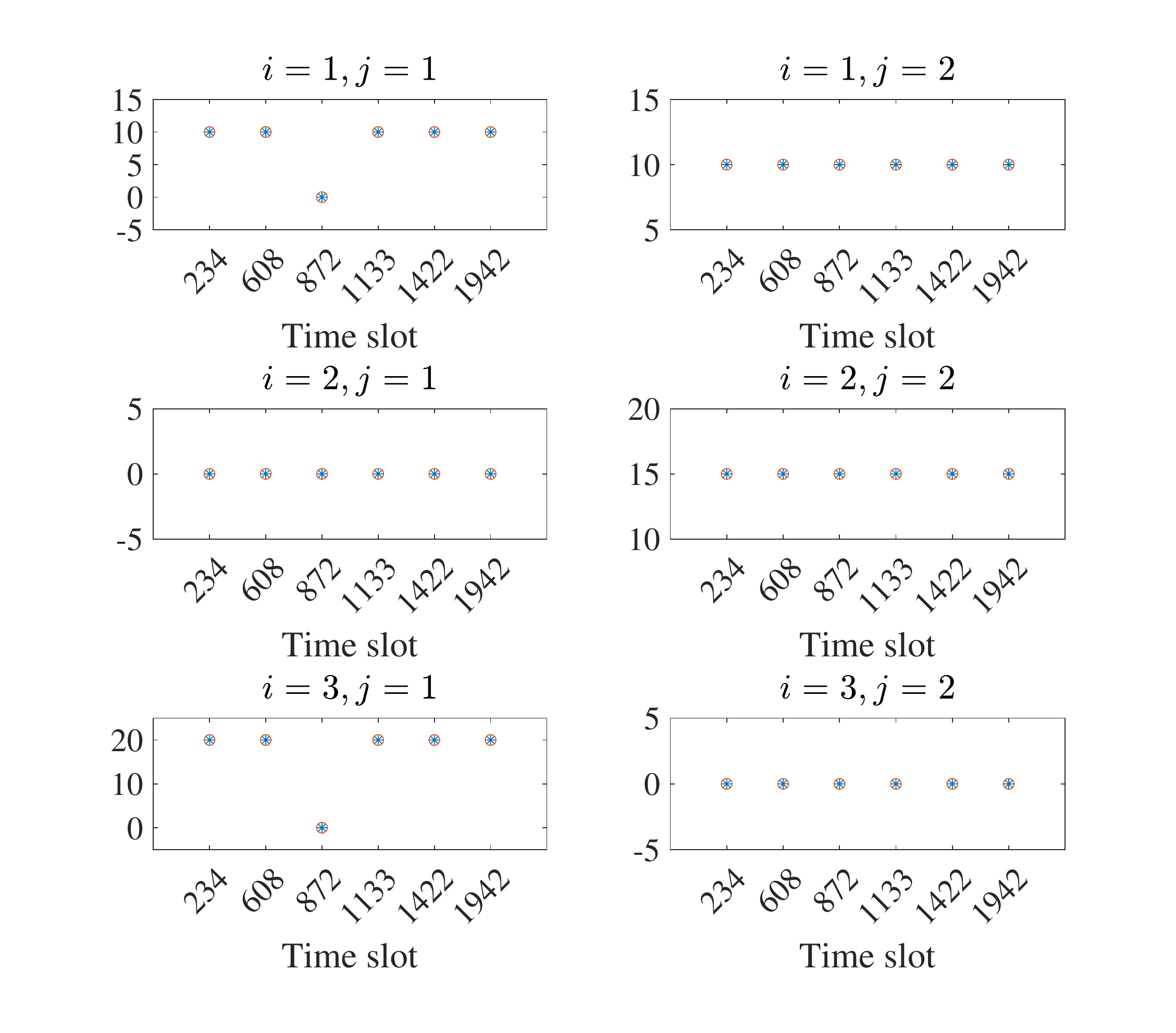}
  \end{minipage}
  \begin{minipage}{0.45\textwidth}
      \includegraphics[width=7.5cm]{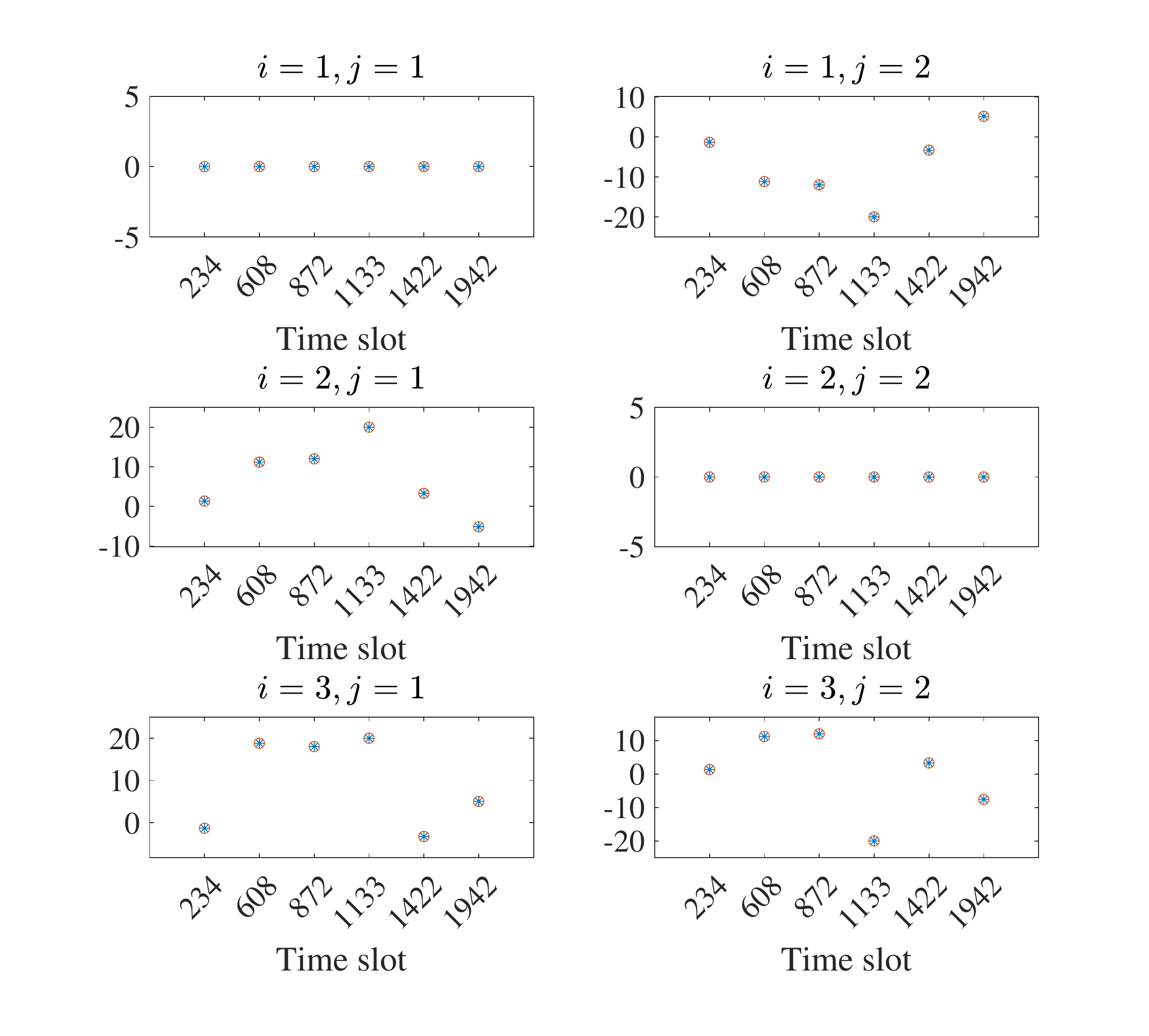}
  \end{minipage}  
  \caption{Left: Comparison of $m_{ij}^{(50)}$ and $(m'_{ij})^{(50)}$ in selected time slots. Blue stars: $m_{ij}^{(50)}$; red circles: $(m'_{ij})^{(50)}$. Right: Comparison of $u_{ik}^{(50)}$ and $(u'_{ik})^{(50)}$ in selected time slots. Blue stars: $u_{ik}^{(50)}$; red circles: $(u'_{ik})^{(50)}$.} \label{fig:m_hat} \label{fig:u_hat}
\end{figure}

In the following, we focus on the effectiveness of the proposed accelerated ADMM algorithm. First, we solve the problem using the traditional ADMM algorithm. To illustrate the computational burden in this test case, the numbers of iterations required to converge are recorded and demonstrated in Figure \ref{fig:histo} as a histogram.
While most of the numbers fall below 10, we also discover that it may take up to 1000 iterations to complete the iteration in extreme cases, probably too time-consuming for the real-time operation of a DC system. 
Note that only 21\% of time slots have gone through more than 50 iterations to converge in Figure \ref{fig:histo}. Therefore, we consider using 50 as the threshold to implement the accelerated algorithm with iteration truncation. We choose 6 examples that take excessive iterations to converge in the traditional algorithm and plot $m_{ij}^{(50)}$ and $(m'_{ij})^{(50)}$, $u_{ik}^{(50)}$ and $(u'_{ik})^{(50)}$ in Figure \ref{fig:m_hat}.
The results show that the discrepancies between $m_{ij}^{(50)}$ and $(m'_{ij})^{(50)}$, $u_{ij}^{(50)}$ and $(u'_{ik})^{(50)}$ are negligible, which means 50 is a proper threshold. The maximum relative errors $\Delta m_{ij}(t)$ and $\Delta u_{ik}(t)$, defined as follows, are less than 0.2\%.
\bsq
\begin{align}
    \Delta m_{ij}(t) & = \left| \frac{m_{ij}^{(N)}(t) - (m'_{ij})^{(N)}(t)}{(m'_{ij})^{(N)}(t)} \right|, \\
    \Delta u_{ik}(t) & = \left| \frac{u_{ik}^{(N)}(t) - (u'_{ik})^{(N)}(t)}{(u'_{ik})^{(N)}(t)} \right|.
\end{align}
\esq

\begin{figure}[htbp]
    \centering
    \begin{minipage}[t]{0.4\textwidth}
        \includegraphics[width=6cm]{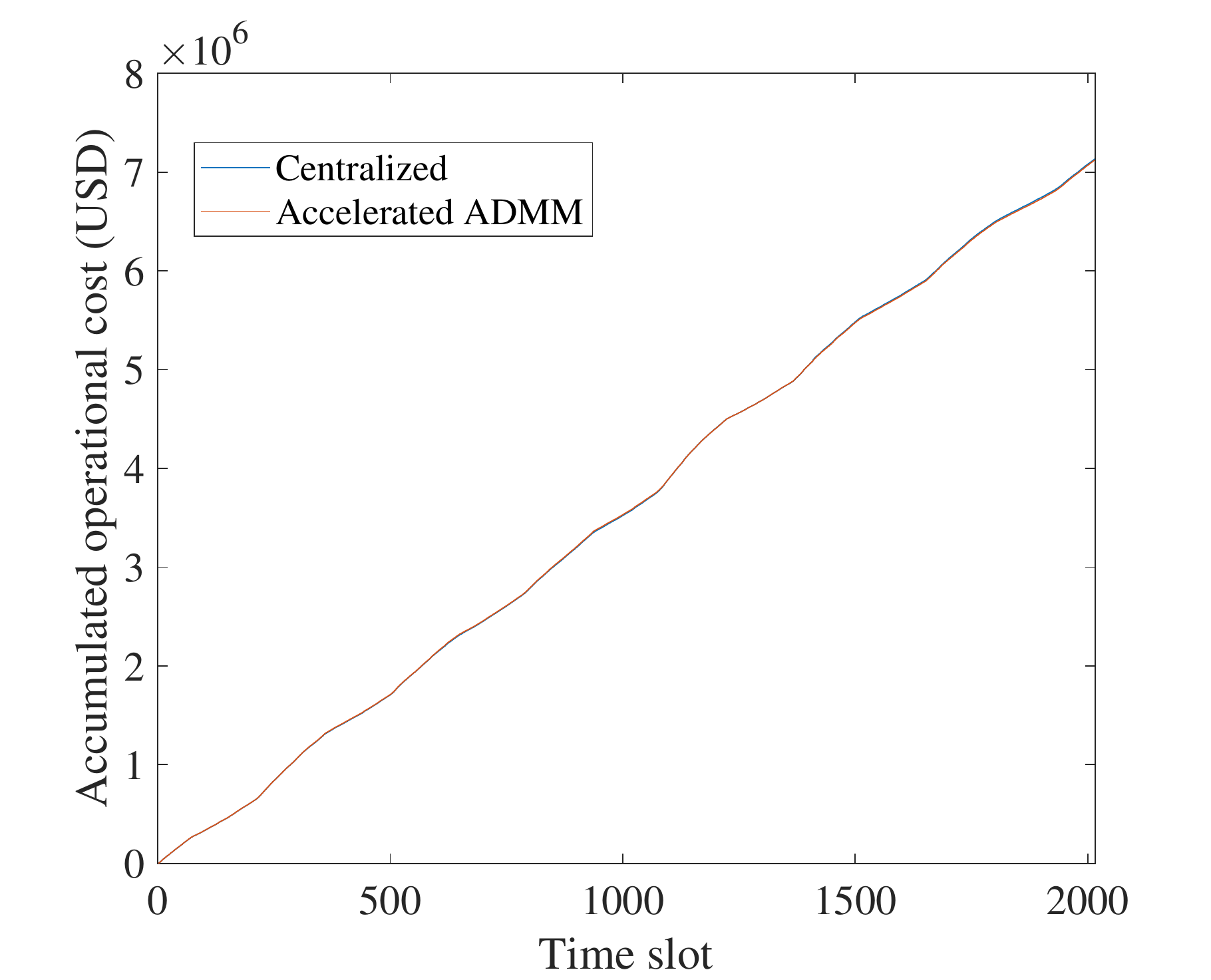}
    \end{minipage}
    \begin{minipage}[t]{0.4\textwidth}
        \includegraphics[width=6cm]{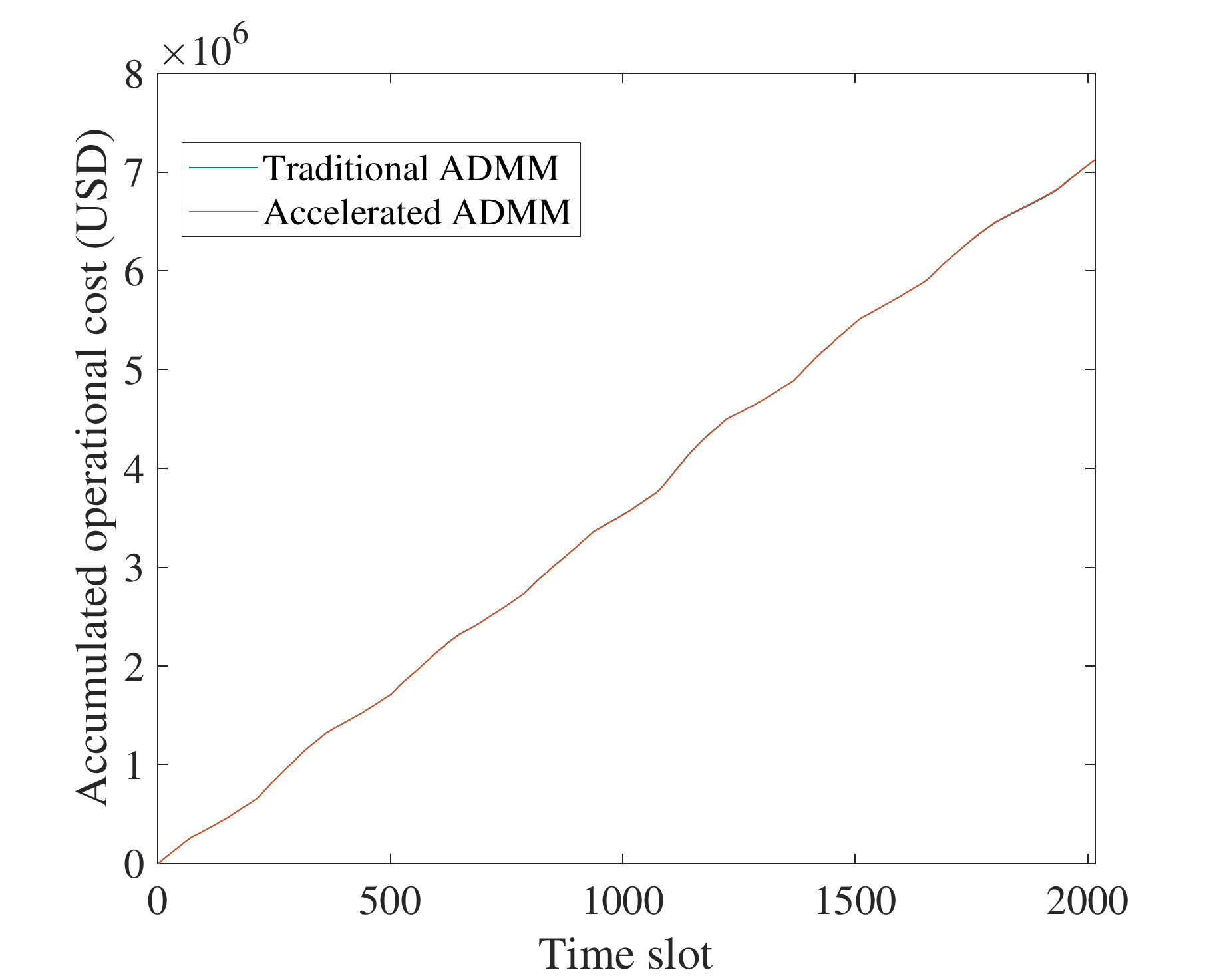}
    \end{minipage}
    \caption{Comparison of cost traces. Left: centralized online algorithm and distributed online algorithm. Right: traditional ADMM algorithm and the proposed accelerated algorithm.} \label{fig:f_hat_1}
\end{figure}

The accumulated overall operation costs obtained by the proposed accelerated ADMM algorithm and its centralized counterpart are compared in Figure \ref{fig:f_hat_1} (left). The two curves coincide with each other, showing the reliability of the proposed distributed algorithm. The results of the traditional ADMM algorithm and the proposed accelerated algorithm with truncation are compared in Figure \ref{fig:f_hat_1} (right). The two traces nearly overlap with each other, which means the truncation does not affect the optimality of the results much.

\begin{figure}[htbp]
  \centering
  \includegraphics[width=0.45\textwidth]{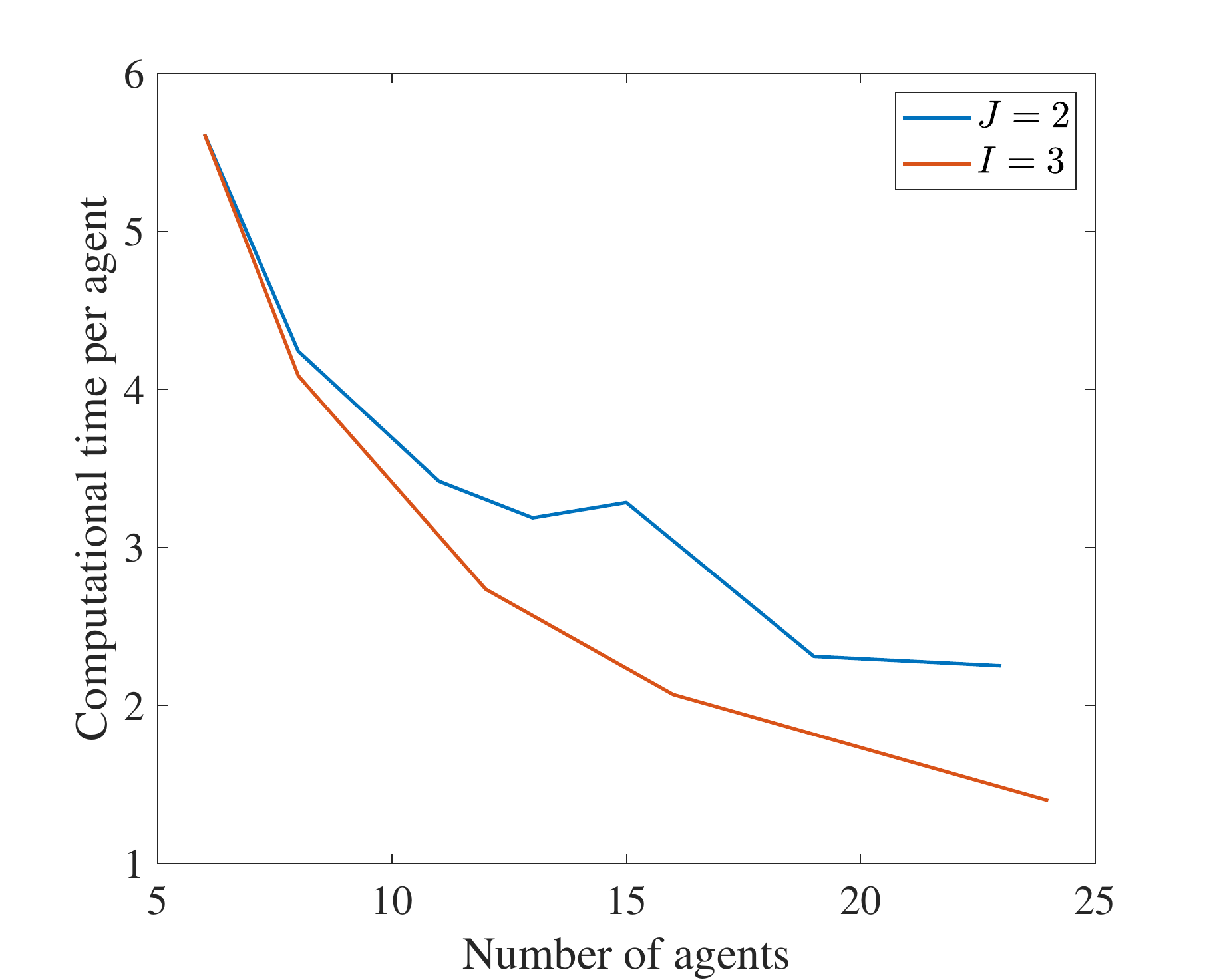} \\
  \caption{Computational time per agent (unit: s).} \label{fig:scal}
\end{figure}

Finally, we evaluate the scalability of the proposed algorithm using larger systems with more front ends and back ends. 
We fix the number of back ends $I$ to 2 and change $J$ from 3 to 20, and fix the number of front ends $J$ to 3 and change $I$ from 2 to 20, respectively. The computational time per agent is shown in Figure \ref{fig:scal}. The time needed is acceptable for real-time DC operation. Moreover, the proposed accelerated ADMM saves more than half (61\%) of the computational time compared to the traditional ADMM.

\section{Conclusion}\label{sec:conclu}
In this paper, a distributed online algorithm is proposed for combined computing workload and energy coordination of data centers. The online algorithm is based on Lyapunov optimization with a novel design of virtual queues. An accelerated ADMM algorithm is developed for fast distributed implementation. Through simulation, we have the following findings:
\begin{enumerate}
    \item The proposed online algorithm is effective and can achieve a nearly offline optimal outcome.
    \item The proposed accelerated ADMM algorithm reduces the computational time compared to the traditional ADMM algorithm.
    \item The proposed algorithm reduces the total costs compared to the greedy algorithm and the one without energy sharing.
\end{enumerate}

Incorporating the heterogeneous features of computing workloads and balancing multiple objectives (total costs, carbon emissions, etc.) will be our future research direction.



\appendix
\makeatletter
\@addtoreset{equation}{section}
\@addtoreset{theorem}{section}
\makeatother
\setcounter{equation}{0}  
\renewcommand{\theequation}{A.\arabic{equation}}
\renewcommand{\thetheorem}{A.\arabic{theorem}}

\section{Proof of Proposition \ref{propo1}} \label{appendix-A}
Before we prove Proposition \ref{propo1}, we give the following lemma.
\begin{lemma}\label{lemma-1}
Suppose $f(z), g(z)$ are differentiable. Suppose $(x^*,y^*,z^*)$ is the optimal solution of the following optimization:
\bsq \label{eq:lemma}
\begin{align}
    \min_{x,y,z}~ & ax-by+f(z) \\
    \mbox{s.t.}~ & x-y=g(z) \\
    ~ & x, y \ge 0, ~\underline{z} \le z \le \overline{z}
\end{align}
\esq
Then, if $ag'(z)+f'(z)\ge 0$ and $bg'(z)+f'(z) \ge 0$, we have $z^*=\underline{z}$.
\end{lemma}

\noindent\emph{Proof of Lemma \ref{lemma-1}}. Problem \eqref{eq:lemma} is equivalent to
\bsq \label{eq:lemma2}
\begin{align}
    \min_{z \in [\underline{z},\overline{z}]} \min_{x,y}~ & ax-by+f(z) \\
    \mbox{s.t.}~ & x-y=g(z) \\
    ~ & x, y \ge 0
\end{align}
\esq

Given $z$, the inner minimization problem is a linear program whose optimal solution lies at one of the vertex of the feasible region. Hence, the optimal solution is either $x^*=0, y^*=-g(z)$ or $x^*=g(z), y^*=0$. If it is the first case, problem \eqref{eq:lemma} can be further turned into
\begin{align}
    \min_{z \in [\underline{z},\overline{z}]}~ bg(z)+f(z)
\end{align}
Therefore, if $bg'(z)+f'(z)\ge 0$, we have $z^*=\underline{z}$. Similarly, if it is the second case, problem \eqref{eq:lemma} can be further turned into
\begin{align}
    \min_{z \in [\underline{z},\overline{z}]}~ ag(z)+f(z)
\end{align}
Therefore, if $ag'(z)+f'(z)\ge 0$, we have $z^*=\underline{z}$. This completes the proof of Lemma \ref{lemma-1}. $\hfill \square$

First, we take the partial derivatives of the objective function in (\ref{obj_online}) with respect to $m_{ij}(t)$, $a_i(t)$, $c_i(t)$, $d_i(t)$, $e_i(t)$, respectively. For $m_{ij}(t)$ and $a_i(t)$ that are decoupled with other variables, we directly calculate their derivatives as follows.
\begin{align}
    \frac{\partial (g(t)+Vf(t))}{\partial m_{ij}(t)}~ & =V \alpha_{ij}-q_j^F(t)+q_i^B(t) + \theta_j-\varphi_i, \\
    \frac{\partial (g(t)+Vf(t))}{\partial a_i(t)}~ & =-V \gamma_j + q_j^F(t)- \theta_j. 
\end{align}

For $c_i(t)$, $d_i(t)$, $e_i(t)$ that are coupled through constraint \eqref{cons:back_end_engy}, eliminating $x_i(t)$ by (\ref{cons:back_end_engy}) and taking the derivatives yields
\begin{align}
    \frac{\partial (g(t)+Vf(t))}{\partial c_i(t)}~ & = V \left(p^\text{b}_i(t) + \beta_i \right) + \eta_\text{c} (b_i(t)-\delta_i-r_iV), \nonumber\\
    \frac{\partial (g(t)+Vf(t))}{\partial d_i(t)} ~ & = V \left(-p^\text{b}_i(t) + \beta_i \right) - \frac{1} {\eta_\text{d}} (b_i(t) - \delta_i-r_iV), \nonumber \\
    \frac{\partial (g(t)+Vf(t))}{\partial e_i(t)} ~ & = Vp^\text{b}_i - q_i^B(t) + \varphi_i.
\end{align}
Similarly, we can eliminate $y_i(t)$ in the same way, yielding
\begin{align}
    \frac{\partial (g(t)+ Vf(t))}{\partial c_i(t)} ~ & = V \left(p^\text{s}_i(t) + \beta_i \right) + \eta_\text{c} (b_i(t) - \delta_i-r_iV),  \nonumber\\
    \frac{\partial (g(t)+ Vf(t))}{\partial d_i(t)} ~ & = V \left(-p^\text{s}_i(t) + \beta_i \right) - \frac{1}{\eta_\text{d}} (b_i(t) - \delta_i-r_iV), \nonumber \\
    \frac{\partial (g(t)+ Vf(t))}{\partial e_i(t)} ~ & = V p^\text{s}_i - q_i^B(t) + \varphi_i. 
\end{align}

After getting the derivatives, we prove the satisfaction of constraints (\ref{cons:qf}), (\ref{cons:qb}), and (\ref{cons:batt_energy}) as follows.

\subsection{Proof of Constraint (\ref{cons:qf})}
\begin{itemize}
    \item When $q_j^F(t) \in \left[0,\sum_{i \in \Omega_j^F} M_{ij}\right]$, we have
    \begin{align}
        \frac{\partial (g(t)+Vf(t))}{\partial m_{ij}(t)} & \ge V\alpha_{ij} - \sum_{i \in \Omega_j^F} M_{ij} + \theta_j - \varphi_i \nonumber \\
        & \ge 0. 
    \end{align}
    The second inequality is derived from the condition that
    \begin{align}
        V \ge \frac{\sum_{i \in \Omega_j^F} M_{ij} - \theta_j + \varphi_i} {\alpha_{ij}}, ~ \forall i \in \Omega_j^F.
    \end{align}
    Thus ${m}_{ji}(t) = 0$ and
    \begin{align}
       0 \le q_j^F(t+1)= ~ & 
 q_j^F(t)+a_j(t)-0 \nonumber\\
 \le ~ &  \sum_{j \in \Omega_j^F} M_{ij} + A_{j,\max} \le Q_j^F,
    \end{align}
which is because of the range of $q_j^F(t)$ in this case and the assumption \textbf{A1}. Therefore, $q_j^F(t+1)$ is still within $[0,Q_j^F]$.
 \item When $q_j^F(t) \in \left[\sum_{i \in \Omega_j^F} M_{ij}, Q_j^F-A_{j,\max}\right]$,
 \bsq
 \begin{align}
     q_j^F(t+1)=~& q_j^F(t)+a_j(t)-\sum_{i \in \Omega_j^F} m_{ij}(t) \nonumber\\
     \ge~ & \sum_{i \in \Omega_j^F} m_{ij}(t)- \sum_{i \in \Omega_j^F} m_{ij}(t) =0, \\
      q_j^F(t+1)=~ &q_j^F(t)+a_j(t)-\sum_{i \in \Omega_j^F} m_{ij}(t) \nonumber\\
      \le ~ &  Q_j^F-A_{j,\max}+A_{j,\max}=Q_j^F.
 \end{align} 
 \esq
    \item When $q_j^F(t) \in \left[Q_j^F-A_{j,\max}, Q_j^F \right]$,
    \begin{align}
        \frac{\partial (g(t) + V f(t))}{\partial a_j(t)} \ge - V \gamma_j + Q_j^F - A_{j,\max} - \theta_j \ge 0.
    \end{align}
    The second inequality is derived from the condition that
    \begin{align}
        V \le \frac{Q_j^F-A_{j,\max}-\theta_j}{\gamma_j}.
    \end{align}
    Thus ${a}_{j}(t) = 0$ and
    \bsq
    \begin{align}
        q_j^F(t+1)=~ & q_j^F(t)+0-\sum_{i \in \Omega_j^F} m_{ij}(t) \le Q_j^F \\
        q_j^F(t+1)=~ & q_j^F(t)+0-\sum_{i \in \Omega_j^F} m_{ij}(t) \nonumber\\
        \ge~ &  Q_j^F-A_{j,\max}-\sum_{i \in \Omega_j^F} M_{ij} \ge 0
    \end{align}
    \esq
which is because of the range of $q_j^F(t)$ in this case and the assumption \textbf{A1}. Therefore, $q_j^F(t+1)$ is still within $[0,Q_j^F]$.
\end{itemize}

\subsection{Proof of Constraint (\ref{cons:qb})} 
\begin{itemize}
    \item When $q_i^B(t) \in \left[0, E_i \right]$, we have
    \bsq
    \begin{align}
        \frac{\partial (g(t)+Vf(t))}{\partial e_i(t)} \ge Vp^\text{b}_i - E_i + \varphi_i \ge 0, \\
        \frac{\partial (g(t)+Vf(t))}{\partial e_i(t)} \ge Vp^\text{s}_i - E_i + \varphi_i \ge 0.
    \end{align}
    \esq
    These can be derived from the condition that
    \begin{align}
         V \ge \frac{E_i-\varphi_i}{p^\text{s}_i} \ge \frac{E_i-\varphi_i}{p^\text{b}_i}.
    \end{align}
    Thus ${e}_{i}(t) = 0$ and
    \begin{align}
        0 \le q_i^B(t+1) = ~& q_i^B(t)+\sum_{j \in \Omega_i^B} m_{ij}(t) \nonumber\\
        \le ~ & E_i + \sum_{j \in \Omega_i^B} M_{ij} \le Q_i^B,
    \end{align}
    which is because of the range of $q_i^B(t)$ in this case, the assumption \textbf{A2}, and Lemma \ref{lemma-1}. Therefore, $q_i^B(t+1)$ is still within $[0, Q_i^B]$.
    \item When $q_i^B(t) \in \left[E_i, Q_i^B-\sum_{j \in \Omega_i^B} M_{ij}\right]$,
    \bsq
    \begin{align}
        q_i^B(t+1) = ~& q_i^B(t)+\sum_{j \in \Omega_i^B} m_{ij}(t)-e_i(t) \nonumber\\
        \ge ~ & E_i-E_i=0, \\
        q_i^B(t+1) = ~& q_i^B(t)+\sum_{j \in \Omega_i^B} m_{ij}(t)-e_i(t) \nonumber\\
        \le ~ & Q_i^B-\sum_{j \in \Omega_i^B} M_{ij}+ \sum_{j \in \Omega_i^B} M_{ij}=Q_i^B.
    \end{align}
    \esq
    \item When $q_i^B(t) \in \left[Q_i^B-\sum_{j \in \Omega_i^B} M_{ij}, Q_i^B\right]$
     \begin{align}
     \frac{\partial (g(t)+Vf(t))}{\partial m_{ij}(t)} \ge ~ & V \alpha_{ij} - Q_j^F + Q_i^B \nonumber \\
     & - \sum_{j \in \Omega_i^B} M_{ij} + \theta_j - \varphi_i \nonumber \\
     \ge ~ & 0
     \end{align}
    The second inequality is derived from
    \begin{align}
        V \ge \frac{Q_j^F-Q_i^B +\sum_{j \in \Omega_i^B} M_{ij} -\theta_j +\varphi_i}{\alpha_{ij}}, \nonumber \\
        \forall j \in \Omega_i^B.
    \end{align}
    Thus ${m}_{ij}(t) = 0$ and
    \bsq
    \begin{align}
        q_i^B(t+1) =~& q_i^B(t)-e_i(t) \le Q_i^B, \\
         q_i^B(t+1) =~ & q_i^B(t)-e_i(t) \nonumber\\
         \ge~& Q_i^B-\sum_{j \in \Omega_i^B} M_{ij}-E_i \ge 0,
    \end{align}
    \esq
    which is because of the range of $q_i^B(t)$ in this case and the assumption \textbf{A2}. 
    
    Therefore, $q_i^B(t+1)$ is still within $[0, Q_i^B]$.
\end{itemize}

\subsection{Proof of Constraint (\ref{cons:batt_energy})}
\begin{itemize}
    \item When $b_i(t) \in \left[\underline{B}_i, \underline{B}_i + \frac{1}{\eta_\text{d}} D_i \right]$, we have
    \bsq
    \begin{align}
        \frac{\partial (g(t) + Vf(t))}{\partial d_i(t)} ~ \ge ~ & V(-p^\text{b}_i(t) + \beta_i) - \frac{1}{\eta_\text{d}}(\underline{B}_i \nonumber \\
        & + \frac{1}{\eta_\text{d} }D_i - r_i V - \delta_i) \nonumber \\
        \ge ~ & 0, \\
        \frac{\partial (g(t) + Vf(t))}{\partial c_i(t)} ~ \ge ~ & V(-p^\text{s}_i(t) + \beta_i) + \frac{1}{\eta_\text{d}}(\underline{B}_i \nonumber \\
        & + \frac{1}{\eta_\text{d} }D_i - r_i V - \delta_i) \nonumber \\
        \ge ~ & 0.
    \end{align}
    \esq
    These can be derived from assumption \textbf{A4} and the conditions that
    \bsq
    \begin{align}
        V \ge \frac{\frac{1}{\eta_\text{d}} ( \underline{B}_i + \frac{1}{\eta_\text{d} }D_i-\delta_i)}{-p^\text{b}_i + \beta_i + \frac{1}{\eta_\text{d}} r_i}, \\
        V \ge \frac{\frac{1}{\eta_\text{d}} ( \underline{B}_i + \frac{1}{\eta_\text{d} }D_i-\delta_i)}{-p^\text{s}_i + \beta_i + \frac{1}{\eta_\text{d}} r_i}.
    \end{align}
    \esq
    Thus $d_i=0$ and
    \begin{align}
       0 \le  b_i(t+1)=~& b_i(t)+\eta_\text{c} c_i(t) \nonumber\\
       \le ~ & \underline{B}_i + \frac{1}{\eta_\text{d}}D_i+\eta_\text{c} C_i \le \overline{B}_i
    \end{align}
    which is because of the range of $b_i(t)$ in this case, the assumption \textbf{A3}, and Lemma \ref{lemma-1}.
    Therefore, $b_i(t+1)$ is still within $[\underline{B}_i,\overline{B}_i]$.
    \item When $b_i(t) \in \left[\underline{B}_i + \frac{1}{\eta_\text{d}} D_i,\overline{B}_i - \eta_\text{c} C_i\right]$,
    \begin{align}
        b_i(t+1)=~ &b_i(t)+\eta_\text{c} c_i(t)-\frac{1}{\eta_\text{d}} d_i(t) \nonumber\\
        \ge ~&  \underline{B}_i + \frac{1}{\eta_\text{d}}D_i - \frac{1}{\eta_\text{d}}D_i = \underline{B}_i \\
        b_i(t+1)=~ &b_i(t)+\eta_\text{c} c_i(t)-\frac{1}{\eta_\text{d}} d_i(t) \nonumber \\
        \le ~ & \overline{B}_i-\eta_\text{c} C_i + \eta_\text{c} C_i = \overline{B}_i
    \end{align}
    \item When $b_i(t) \in \left[\overline{B}_i - \eta_\text{c} C_i, \overline{B}_i \right]$, we have
    \bsq
    \begin{align}
        \frac{\partial (g(t)+Vf(t))}{\partial c_i(t)} ~ \ge ~ & V(p^\text{b}_i(t)+\beta_i) \nonumber \\
        & + \eta_\text{c}(\overline{B}_i - \eta_\text{c} C_i - \delta_i-\gamma_iV) \nonumber \\
        \ge ~ & 0 \\
        \frac{\partial (g(t)+Vf(t))}{\partial c_i(t)} ~ \ge ~ & V(p^\text{s}_i(t)+\beta_i) \nonumber \\
        & + \eta_\text{c}(\overline{B}_i - \eta_\text{c} C_i - \delta_i-\gamma_iV) \nonumber \\
        \ge ~ & 0
    \end{align}
    \esq
    These can be derived from assumption \textbf{A4} and the conditions that
    \bsq
    \begin{align}
         V \ge \frac{\eta_\text{c}(\delta_i + \eta_\text{c} C_i - \overline{B}_i)}{p^\text{b}_i + \beta_i - \eta_\text{c} r_i} \\
         V \ge \frac{\eta_\text{c}(\delta_i + \eta_\text{c} C_i - \overline{B}_i)}{p^\text{s}_i + \beta_i - \eta_\text{c} r_i}
    \end{align}
    \esq
    Thus $\hat{c}_i = 0$ and
    \begin{align}
        b_i(t+1)=~&b_i(t)-\frac{1}{\eta_\text{d}} d_i(t) \le \overline{B}_i, \nonumber\\
        b_i(t+1)=~ & b_i(t)-\frac{1}{\eta_\text{d}} d_i(t) \ge \overline{B}_i-\eta_\text{c} C_i - \frac{1}{\eta_\text{d}} D_i \nonumber\\
        \ge~ & \underline{B}_i.
    \end{align}
    which is because of the range of $b_i(t)$ in this case, the assumption \textbf{A3}, and Lemma \ref{lemma-1}.
    
    Therefore, $b_i(t+1)$ is still within $[\underline{B}_i,\overline{B}_i]$.
\end{itemize}

This completes the proof. $\hfill \square$

\renewcommand{\theequation}{B.\arabic{equation}}
\renewcommand{\thetheorem}{B.\arabic{theorem}}

\section{Proof of Proposition \ref{propo2}}
\label{appendix-B}
According to (\ref{ineq3}), we have
\begin{align}\label{ineq4}
    & \quad ~ \Delta(\boldsymbol{\Theta}(t)) + V\mathbb{E} \left[\hat{f}(t) | \boldsymbol{\Theta}(t) \right] \nonumber \\
    &\le \sum_j \mathbb{E}[h_j^F(t)|\boldsymbol{\Theta}(t) ] \mathbb{E}\left[\left(\hat{a}_j(t) - \sum_{i \in \Omega_j^F} \hat{m}_{ji}(t)\right) |~ \boldsymbol{\Theta}(t) \right] \nonumber \\
    & \quad + \sum_i \mathbb{E}[h_i^B(t)|\boldsymbol{\Theta}(t)] \mathbb{E}\left[\left(\sum_{j \in \Omega_i^B} \hat{m}_{ji}(t) - \hat{e}_i(t)\right) |~ \boldsymbol{\Theta}(t) \right] \nonumber \\
    &\quad+ \sum_i \mathbb{E}[l_i(t)|\boldsymbol{\Theta}(t) ] \mathbb{E}\left[\left(\eta_\text{c} \hat{c}_i(t) - \frac{1}{\eta_\text{d}} \hat{d}_i(t)\right) |~ \boldsymbol{\Theta}(t)\right] \nonumber \\
    &\quad+ V\mathbb{E} \left[\hat{f}(t) | \boldsymbol{\Theta}(t) \right] + N_1 + N_2 + N_3 \nonumber \\
    &\le \sum_j \mathbb{E}[h_j^F(t)|\boldsymbol{\Theta}(t) ] \mathbb{E}\left[\left(a_j^*(t) - \sum_{i \in \Omega_j^F} m_{ij}^*(t)\right) \right] \nonumber \\
    &\quad+ \sum_i \mathbb{E}[h_i^B(t)|\boldsymbol{\Theta}(t) ] \mathbb{E}\left[\left(\sum_{j \in \Omega_i^B} m_{ij}^*(t) - e_i^*(t)\right) \right] \nonumber \\
    &\quad+ \sum_i \mathbb{E}[l_i(t)|\boldsymbol{\Theta}(t) ] \mathbb{E}\left[\left(\eta_\text{c} c_i^*(t) - \frac{1}{\eta_\text{d}} d_i^*(t)\right) \right] \nonumber \\
    &\quad+ V\mathbb{E} \left[f^*(t)\right] + N_1 + N_2 + N_3
\end{align}

According to the strong law of large numbers, 
\begin{align}
    & \mathbb{E}\left[\left(a_j^*(t) - \sum_{i \in \Omega_j^F} m_{ij}^*(t)\right) \right] \nonumber \\
    = & \lim_{T \to \infty} \frac{1}{T} \sum_{t = 1}^T \left(a_j^*(t) - \sum_{i \in \Omega_j^F} m_{ij}^*(t)\right) \nonumber \\
    = & 0 \label{eq2a} \\
    & \mathbb{E}\left[\left(\sum_{j \in \Omega_i^B} m_{ij}^*(t) - e_i^*(t)\right) \right] \nonumber \\
    = & \lim_{T \to \infty} \frac{1}{T} \sum_{t = 1}^T \left(\sum_{j \in \Omega_i^B} m_{ij}^*(t) - e_i^*(t)\right) \nonumber \\
    = & 0 \label{eq2b} \\
    & \mathbb{E}\left[\left(\eta_\text{c} c_i^*(t) - \frac{1}{\eta_\text{d}} d_i^*(t)\right) \right] \nonumber \\
    = & \lim_{T \to \infty} \frac{1}{T} \sum_{t = 1}^T \left(\eta_\text{c} c_i^*(t) - \frac{1}{\eta_\text{d}} d_i^*(t)\right) \nonumber \\
    = & 0 \label{eq2c}
\end{align}
The second equality of (\ref{eq2a}) -- (\ref{eq2c}) can be derived from (\ref{vir_que1}) -- (\ref{vir_que3}), hence the right-hand side of (\ref{ineq4}) can be reduced. By summing the new inequality over $t = 1, 2, ..., T$, we have
\begin{align} \label{ineq2}
    & L(\boldsymbol{\Theta}(t + 1)) - L(\boldsymbol{\Theta}(1)) + V \sum_{t = 1}^T \mathbb{E} \left[\hat{f}(t)\right] \nonumber \\
    \le & V \sum_{t = 1}^T \mathbb{E} \left[f^*(t)\right] + T\left(N_1 + N_2 + N_3\right)
\end{align}
We divide both sides of (\ref{ineq2}) by $VT$ and take the limit with $T$ going to infinity, yielding
\begin{align}
    & \underbrace{\lim_{T \to \infty} \frac{1}{T} \sum_{t = 1}^T \mathbb{E} \left[\hat{f}(t)\right]}_{\hat F} - \underbrace{\lim_{T \to \infty} \frac{1}{T} \sum_{t = 1}^T \mathbb{E} \left[f^*(t)\right]}_{F^*} \nonumber \\
    \le & \frac{1}{V}\left(N_1 + N_2 + N_3\right)
\end{align}
This completes the proof. $\hfill \square$

\bibliographystyle{elsarticle-num} 
\bibliography{refs}

\begin{thebibliography}{10}
\expandafter\ifx\csname url\endcsname\relax
  \def\url#1{\texttt{#1}}\fi
\expandafter\ifx\csname urlprefix\endcsname\relax\def\urlprefix{URL }\fi
\expandafter\ifx\csname href\endcsname\relax
  \def\href#1#2{#2} \def\path#1{#1}\fi

\bibitem{miller2020sustainability}
R.~Miller, The sustainability imperative: Green data centers and our cloudy
  future, Data Center Frontier, Tech. Rep (2020).

\bibitem{regional}
R.~Miller, Regional data center clusters power amazon’s cloud[Online].
  Available:
  \url{https://www.datacenterfrontier.com/featured/article/11431479/regional-data-center-clusters-power-amazon8217s-cloud}
  (2015).

\bibitem{WANG2022107926}
H.~Wang, Q.~Wang, Y.~Tang, Y.~Ye, Spatial load migration in a power system:
  Concept, potential and prospects, International Journal of Electrical Power
  \& Energy Systems 140 (2022) 107926.
\newblock \href {https://doi.org/https://doi.org/10.1016/j.ijepes.2021.107926}
  {\path{doi:https://doi.org/10.1016/j.ijepes.2021.107926}}.

\bibitem{ALOBAIDI2021107231}
A.~H. Alobaidi, M.~Khodayar, A.~Vafamehr, H.~Gangammanavar, M.~E. Khodayar,
  Stochastic expansion planning of battery energy storage for the
  interconnected distribution and data networks, International Journal of
  Electrical Power \& Energy Systems 133 (2021) 107231.
\newblock \href {https://doi.org/https://doi.org/10.1016/j.ijepes.2021.107231}
  {\path{doi:https://doi.org/10.1016/j.ijepes.2021.107231}}.

\bibitem{Ding}
Z.~Ding, L.~Xie, Y.~Lu, P.~Wang, S.~Xia, Emission-aware stochastic resource
  planning scheme for data center microgrid considering batch workload
  scheduling and risk management, IEEE Transactions on Industry Applications
  54~(6) (2018) 5599--5608.
\newblock \href {https://doi.org/10.1109/TIA.2018.2851516}
  {\path{doi:10.1109/TIA.2018.2851516}}.

\bibitem{WangP}
P.~Wang, Y.~Cao, Z.~Ding, H.~Tang, X.~Wang, M.~Cheng, Stochastic programming
  for cost optimization in geographically distributed internet data centers,
  CSEE Journal of Power and Energy Systems 8~(4) (2022) 1215--1232.
\newblock \href {https://doi.org/10.17775/CSEEJPES.2020.02930}
  {\path{doi:10.17775/CSEEJPES.2020.02930}}.

\bibitem{NIU2021106358}
T.~Niu, B.~Hu, K.~Xie, C.~Pan, H.~Jin, C.~Li, Spacial coordination between data
  centers and power system considering uncertainties of both source and load
  sides, International Journal of Electrical Power \& Energy Systems 124 (2021)
  106358.
\newblock \href {https://doi.org/https://doi.org/10.1016/j.ijepes.2020.106358}
  {\path{doi:https://doi.org/10.1016/j.ijepes.2020.106358}}.

\bibitem{Chen}
T.~Chen, Y.~Zhang, X.~Wang, G.~B. Giannakis, Robust workload and energy
  management for sustainable data centers, IEEE Journal on Selected Areas in
  Communications 34~(3) (2016) 651--664.
\newblock \href {https://doi.org/10.1109/JSAC.2016.2525618}
  {\path{doi:10.1109/JSAC.2016.2525618}}.

\bibitem{Jawad}
M.~Jawad, M.~B. Qureshi, M.~U.~S. Khan, S.~M. Ali, A.~Mehmood, B.~Khan,
  X.~Wang, S.~U. Khan, A robust optimization technique for energy cost
  minimization of cloud data centers, IEEE Transactions on Cloud Computing
  9~(2) (2021) 447--460.
\newblock \href {https://doi.org/10.1109/TCC.2018.2879948}
  {\path{doi:10.1109/TCC.2018.2879948}}.

\bibitem{Fang}
Q.~Fang, J.~Wang, Q.~Gong, Qos-driven power management of data centers via
  model predictive control, IEEE Transactions on Automation Science and
  Engineering 13~(4) (2016) 1557--1566.
\newblock \href {https://doi.org/10.1109/TASE.2016.2582501}
  {\path{doi:10.1109/TASE.2016.2582501}}.

\bibitem{Mansouri}
Y.~Mansouri, A.~N. Toosi, R.~Buyya, Cost optimization for dynamic replication
  and migration of data in cloud data centers, IEEE Transactions on Cloud
  Computing 7~(3) (2019) 705--718.
\newblock \href {https://doi.org/10.1109/TCC.2017.2659728}
  {\path{doi:10.1109/TCC.2017.2659728}}.

\bibitem{cupelli2018data}
L.~Cupelli, T.~Sch{\"u}tz, P.~Jahangiri, M.~Fuchs, A.~Monti, D.~M{\"u}ller,
  Data center control strategy for participation in demand response programs,
  IEEE Transactions on Industrial Informatics 14~(11) (2018) 5087--5099.

\bibitem{Guo}
Y.~Guo, Y.~Gong, Y.~Fang, P.~P. Khargonekar, X.~Geng, Energy and network aware
  workload management for sustainable data centers with thermal storage, IEEE
  Transactions on Parallel and Distributed Systems 25~(8) (2014) 2030--2042.
\newblock \href {https://doi.org/10.1109/TPDS.2013.278}
  {\path{doi:10.1109/TPDS.2013.278}}.

\bibitem{WANG2021106451}
G.~Wang, X.~Yang, W.~Cai, Y.~Zhang, Event-triggered online energy flow control
  strategy for regional integrated energy system using lyapunov optimization,
  International Journal of Electrical Power \& Energy Systems 125 (2021)
  106451.
\newblock \href {https://doi.org/https://doi.org/10.1016/j.ijepes.2020.106451}
  {\path{doi:https://doi.org/10.1016/j.ijepes.2020.106451}}.

\bibitem{Karimiafshar}
A.~Karimiafshar, M.~R. Hashemi, M.~R. Heidarpour, A.~N. Toosi, Effective
  utilization of renewable energy sources in fog computing environment via
  frequency and modulation level scaling, IEEE Internet of Things Journal
  7~(11) (2020) 10912--10921.
\newblock \href {https://doi.org/10.1109/JIOT.2020.2993276}
  {\path{doi:10.1109/JIOT.2020.2993276}}.

\bibitem{Sun}
J.~Sun, S.~Chen, P.~You, Q.~Yang, Z.~Yang, Battery-assisted online operation of
  distributed data centers with uncertain workload and electricity prices, IEEE
  Transactions on Cloud Computing (2021) 1--1\href
  {https://doi.org/10.1109/TCC.2021.3132174}
  {\path{doi:10.1109/TCC.2021.3132174}}.

\bibitem{SONG2022108289}
M.~Song, Y.~Cai, C.~Gao, T.~Chen, Y.~Yao, H.~Ming, Transactive energy in power
  distribution systems: Paving the path towards cyber-physical-social system,
  International Journal of Electrical Power \& Energy Systems 142 (2022)
  108289.
\newblock \href {https://doi.org/https://doi.org/10.1016/j.ijepes.2022.108289}
  {\path{doi:https://doi.org/10.1016/j.ijepes.2022.108289}}.

\bibitem{Yu2}
L.~Yu, T.~Jiang, Y.~Zou, Distributed real-time energy management in data center
  microgrids, IEEE Transactions on Smart Grid 9~(4) (2018) 3748--3762.
\newblock \href {https://doi.org/10.1109/TSG.2016.2640453}
  {\path{doi:10.1109/TSG.2016.2640453}}.

\bibitem{yan2022distributed}
D.~Yan, C.~Ma, Y.~Chen, Distributed coordination of charging stations
  considering aggregate ev power flexibility, IEEE Transactions on Sustainable
  Energy 14~(1) (2022) 356--370.

\bibitem{xu2020peer}
D.~Xu, B.~Zhou, N.~Liu, Q.~Wu, N.~Voropai, C.~Li, E.~Barakhtenko, Peer-to-peer
  multienergy and communication resource trading for interconnected microgrids,
  IEEE Transactions on Industrial Informatics 17~(4) (2020) 2522--2533.

\bibitem{ZHANG2020106094}
R.~Zhang, K.~Yan, G.~Li, T.~Jiang, X.~Li, H.~Chen, Privacy-preserving
  decentralized power system economic dispatch considering carbon capture power
  plants and carbon emission trading scheme via over-relaxed admm,
  International Journal of Electrical Power \& Energy Systems 121 (2020)
  106094.
\newblock \href {https://doi.org/https://doi.org/10.1016/j.ijepes.2020.106094}
  {\path{doi:https://doi.org/10.1016/j.ijepes.2020.106094}}.

\bibitem{RAJAEI2021107126}
A.~Rajaei, S.~Fattaheian-Dehkordi, M.~Fotuhi-Firuzabad, M.~Moeini-Aghtaie,
  Decentralized transactive energy management of multi-microgrid distribution
  systems based on admm, International Journal of Electrical Power \& Energy
  Systems 132 (2021) 107126.
\newblock \href {https://doi.org/https://doi.org/10.1016/j.ijepes.2021.107126}
  {\path{doi:https://doi.org/10.1016/j.ijepes.2021.107126}}.

\bibitem{ullah2021peer}
M.~H. Ullah, J.-D. Park, Peer-to-peer energy trading in transactive markets
  considering physical network constraints, IEEE Transactions on Smart Grid
  12~(4) (2021) 3390--3403.

\bibitem{pjm}
{PJM-Data Miner 2}, Real-time five minute lmps[Online]. Available:
  \url{https://dataminer2.pjm.com} (2022).

\bibitem{google}
{John Wilkes and Charles Reiss}, Google clusterdata 2011 traces[Online].
  Available:
  \url{https://github.com/google/cluster-data/blob/master/ClusterData2011_2.md}
  (2021).

\bibitem{solardata}
{National Renewable Energy Laboratory}, National solar radiation database
  data[Online]. Available:
  \url{https://nsrdb.nrel.gov/data-sets/archives.html/} (2018).

\end{thebibliography}
\end{document}